\documentclass[12pt,reqno]{amsart}
\usepackage{amsmath,amsfonts,amsthm,amsopn,amssymb,mathtools,mathrsfs,enumerate,color}
\usepackage{amsmath}
\usepackage{comment}

\usepackage{cite,marginnote}
\usepackage{mathtools} 

\usepackage{color,graphicx,enumerate}
\usepackage[colorlinks=true,urlcolor=blue,
citecolor=red,linkcolor=blue,linktocpage,pdfpagelabels,
bookmarksnumbered,bookmarksopen]{hyperref}
\usepackage[english]{babel}

\usepackage[left=2.9cm,right=2.9cm,top=2.8cm,bottom=2.8cm]{geometry}

\numberwithin{equation}{section}

\makeindex
\newtheorem{theorem}{Theorem}[section]

\newtheorem{remark}{Remark}[section]

\newtheorem{proposition}{Proposition}[section]
\newtheorem{lemma}{Lemma}[section]
\newtheorem{iteration lemma}{iteration Lemma}[section]

\newcommand{\R}{\mathbb R}

\title[positive solutions of the Choquard equation]
{Symmetry, uniqueness and nondegeneracy of positive solutions of the Choquard equation}

\author[S.~Fan]{Song Fan}
\author[G-D.~Li]{Gui-Dong Li}
\author[J.~J.~Zhang]{Jianjun Zhang}

\address[S.~Fan]{\newline\indent School  of  Mathematics  and  Statistics
\newline\indent
Guizhou University
\newline\indent
Guiyang, 550025, Guizhou, PR China}
\email{\href{mailto:doraemonsong77@gmail.com}{doraemonsong77@gmail.com}}

\address[G-D.~Li]{\newline\indent School  of  Mathematics  and  Statistics
\newline\indent
Guizhou University
\newline\indent
Guiyang, 550025, Guizhou, PR China}
\email{\href{mailto:bestdong123@163.com}{bestdong123@163.com}}

\address[J.~J.~Zhang]{\newline\indent College of Mathematics and Statistics
\newline\indent
Chongqing Jiaotong University
\newline\indent
Xuefu, Nan'an, 400074, Chongqing, PR China}
\email{\href{mailto:zhangjianjun09@tsinghua.org.cn}{zhangjianjun09@tsinghua.org.cn}}

\begin{document}
\begin{abstract}
  We study positive solutions \(Q\in H^1(\mathbb R^N)\) of
\[
    -\Delta Q+Q=(I_\alpha*Q^p)Q^{p-1},
    \]
 where \(I_\alpha\) is the Riesz potential of order \(\alpha\), \(N\ge1\),
   \(0<\alpha<N\), \(p>1+\alpha/N\), and \((N-2)p<N+\alpha\). We prove that every such
   solution is radially symmetric and strictly decreasing about some
   point. For \(p\ge2\), we establish uniqueness up to translations and show
 that the real linearized operator has Morse index one and kernel spanned by
   \(\partial_1Q,\ldots,\partial_NQ\).

  The symmetry proof uses the method of moving planes. Uniqueness follows from nondegeneracy and
   compactness through a continuation argument.

   \vskip0.23in


\end{abstract}

\maketitle

\section{Introduction}
 We study positive solutions of the Choquard equation
    \begin{equation}\label{eq:equation}
  -\Delta Q+Q=(I_\alpha*Q^p)Q^{p-1} \quad\text{in }\R^N,\qquad Q>0,\quad Q\in H^1(\R^N).
  \end{equation}
    Here \(I_\alpha\) is the Riesz potential, normalized by
  \[
  I_\alpha(x)=c_{N,\alpha}|x|^{\alpha-N},\qquad c_{N,\alpha} =\frac{\Gamma((N-\alpha)/2)} {2^\alpha\pi^{N/2}\Gamma(\alpha/2)}.
  \]
  We consider the parameter range
  \begin{equation}\label{eq:existence-range}
    N\ge1,\qquad 0<\alpha<N,\qquad p>1+\frac\alpha N,\qquad (N-2)p<N+\alpha.
\end{equation}
  Under~\eqref{eq:existence-range}, equation \eqref{eq:equation} is well defined on
  \(H^1(\R^N)\) by the Hardy--Littlewood--Sobolev (HLS) and Sobolev
    inequalities~\cite{Lieb83,MV13}.

  The case \((N,\alpha,p)=(3,2,2)\), for which
  \(I_2(x)=(4\pi|x|)^{-1}\), arises in Pekar's polaron
  model~\cite{Pekar} and in Choquard's approximation to
  Hartree--Fock theory~\cite{Lieb76}. Its variational theory was
  developed by Lieb~\cite{Lieb76} and Lions~\cite{Lions80}. More
  generally, equation \eqref{eq:equation} is the standing-wave equation for the
  focusing nonlinear Hartree equation. We refer to~\cite{MV17} for its
  history and physical interpretations.

 A ground state is a nonzero solution whose action is minimal among all nonzero solutions. Moroz and
 Van Schaftingen~\cite{MV13} proved the existence of positive ground states
 throughout~\eqref{eq:existence-range}, together with their regularity and radial symmetry up to translations. For
 sufficiently regular nonzero solutions, the Pohozaev identity also gives the necessity of both strict inequalities
 in~\eqref{eq:existence-range}.

For \(p\ge2\), we study uniqueness up to translations among all
positive \(H^1(\R^N)\) solutions, together with nondegeneracy of
each such solution. The real linearization at a positive solution
\(Q\) is
  \begin{equation}\label{eq:linearization}
  L_Qh=-\Delta h+h-(p-1)(I_\alpha*Q^p)Q^{p-2}h -pQ^{p-1}\bigl(I_\alpha*(Q^{p-1}h)\bigr).
\end{equation}
  Translation invariance gives \(L_Q\partial_jQ=0\). Nondegeneracy means that 
  \[\ker L_Q=\operatorname{span}\{\partial_1Q,\ldots,\partial_NQ\}. 
  \]
  The Morse index of \(L_Q\) is the
  number of its negative eigenvalues, counted with multiplicity.

  For \((N,\alpha,p)=(3,2,2)\), Lieb~\cite{Lieb76} proved uniqueness, up to translations, of the positive ground state.
  Lenzmann~\cite{Lenzmann09} and Wei and
 Winter~\cite{WW09} proved nondegeneracy of the positive
 radial ground state.  Ma and Zhao~\cite{MZ10} classified all positive solutions.
 For \(\alpha=2\)
 and \(p=2\), the ODE analysis of Choquard, Stubbe and Vuffray~\cite{CSV08} also yields
 uniqueness of positive radial \(H^1(\R^N)\) solutions in dimensions \(N\in\{3,4,5\}\).
 Krieger, Lenzmann and Rapha\"el~\cite{KLR09} established uniqueness of positive solution and radial nondegeneracy in
 dimension four. Arora and Roudenko~\cite{AR22} presented a proof of radial uniqueness of positive solution for
 \(N\in\{3,4,5\}\), and Chen~\cite{Chen21} proved ground-state nondegeneracy in \(N\in\{3,4,5\}\).

 Perturbative results give uniqueness and nondegeneracy near
  several parameter values. Xiang~\cite{Xiang16} proved both
  properties for ground states when \(N=3\), \(\alpha=2\), and
 \(2<p<2+\delta\), for some \(\delta>0\). Luo, Zhang and
  Xu~\cite{LZX24} obtained the same conclusions for ground states
  when \(N=3\), \(p=2\), and \(2-\varepsilon<\alpha<2+\varepsilon\), for
  some \(\varepsilon\in(0,1)\). Seok~\cite{Seok19} treated the limits
  \(\alpha\downarrow0\) and \(\alpha\uparrow N\). For each fixed pair
  \((N,p)\) in the corresponding range, ground-state uniqueness and
  nondegeneracy hold for some \(\varepsilon=\varepsilon(N,p)>0\) when
\[
   \begin{aligned}
 &N=3,\qquad 2\le p<3, &&0<\alpha<\varepsilon,\\
&N\ge3,\qquad 2<p<\frac{2N}{N-2}, &&N-\varepsilon<\alpha<N.
    \end{aligned}
 \]
  Li~\cite{Li22} proved uniqueness and nondegeneracy for all
   positive solutions near the Newtonian quadratic case: for each
   \(N\in\{3,4,5\}\), there exists \(\delta=\delta(N)>0\) such that both
   conclusions hold for
  \[
|\alpha-2|\le\delta,\qquad 2\le p\le2+\delta.
  \]

\subsection*{Main results}

 We first establish radial symmetry throughout the existence
   range~\eqref{eq:existence-range}.

\begin{theorem}[Radial symmetry]
    \label{thm:radial-symmetry}
 Let \( N\ge1\), \(0<\alpha<N\),  \(p> 1+\frac{\alpha}{N}\), \((N-2)p<N+\alpha\). Let \(u\in H^1(\R^N)\) be a positive  solution
      of~\eqref{eq:equation}. Then there exist \(x_0\in\R^N\) and \(U\in C^2([0,\infty))\) such that
 \[
 u(x)=U(|x-x_0|) \quad\text{for all }x\in\R^N, \qquad U'(r)<0 \quad\text{for all }r>0.
   \]
\end{theorem}
Moroz and Van Schaftingen~\cite{MV13} proved that every positive ground state
in \eqref{eq:existence-range} is radial and radially nonincreasing up to
translation. Their proof uses polarization and the minimizing property of a
ground state: polarization preserves its \(L^2\) and gradient norms and increases
the Riesz interaction; equality in the resulting variational inequality yields
symmetry.

 That argument does not establish symmetry for a positive solution which is
 not known to be a ground state. For \(p\ge2\), symmetry of all positive solutions
 follows from \cite{BSR26}. For \(p<2\) within \eqref{eq:existence-range}, we
 prove symmetry of all positive \(H^1(\R^N)\) solutions by moving
 planes in Section~\ref{sec:symmetry-section}; see~\cite{BN91,CLO05,CLL17} for
 this method in local and nonlocal equations.

 Our results on the Morse index, nondegeneracy, and uniqueness of
 positive solutions concern the parameter range
  \begin{equation}\label{eq:range}
   N\ge1,\qquad 0<\alpha<N,\qquad p\ge2,\qquad (N-2)p<N+\alpha.
   \end{equation}

\begin{theorem}[Morse index and nondegeneracy]
\label{thm:index-one-nondegeneracy}
Let \( N\ge1\), \(0<\alpha<N\),  \(p\ge2\), \((N-2)p<N+\alpha\). Then every positive  solution \(Q\in H^1(\R^N)\) of \eqref{eq:equation}
   satisfies
 \[
    \operatorname{ind}(L_Q)=1, \qquad \ker L_Q
    =\operatorname{span}\{\partial_1Q,\ldots,\partial_NQ\},
 \]
  where \(L_Q\) is the real linearized operator defined in \eqref{eq:linearization}.
  Thus every positive solution has Morse index one and is nondegenerate modulo translations.
\end{theorem}

\begin{theorem}[Uniqueness of positive solutions]
  \label{thm:unique}
Let \( N\ge1\), \(0<\alpha<N\),  \(p\ge2\), \((N-2)p<N+\alpha\). Then the positive radial \(H^1(\R^N)\) solution of \eqref{eq:equation}
 is unique; we denote it by \(Q\). Moreover, every positive \(H^1(\R^N)\) solution of
  \eqref{eq:equation} is of the form
 \[
 Q(\cdot-x_0),\qquad x_0\in\R^N.
  \]
\end{theorem}

  The proof of Theorem~\ref{thm:unique} combines
  Theorem~\ref{thm:index-one-nondegeneracy} with a
  continuation argument. By Theorem~\ref{thm:radial-symmetry}, each
  positive solution can be translated to a radial decreasing one. For
  fixed \(N\), Li~\cite[Theorem~1.5]{Li22} then gives uniform
  \(W^{2,r}(\R^N)\) and exponential-decay bounds as \((\alpha,p)\) varies
  over compact subsets of the parameter set in \eqref{eq:range}. These
  estimates provide the compactness needed to continue the solutions in the
  parameters. This use of nondegeneracy and continuation is closely related
  to the proofs of ground-state uniqueness in~\cite{FL13,FLS16}; see
  also~\cite{FV14,ADG22,DP25} for related arguments.

\subsection*{Proof strategy}

  \emph{Radial symmetry.} We use the method of moving planes. For \(p\ge2\), the result follows from B\"oer, dos
  Santos and Ramos \cite[Theorem~1.1]{BSR26}. For \(1+\alpha/N<p<2\), the comparison requires an estimate that
  accounts for the negative exponent \(p-2\).

 Write \(x=(x_1,x')\) and set
 \[
  \begin{aligned}
  &\Sigma_\lambda=\{x_1<\lambda\},\qquad x^\lambda=(2\lambda-x_1,x'),\qquad
  u_\lambda(x)=u(x^\lambda),\\
 &E_\lambda=\{x\in\Sigma_\lambda:u(x)>u_\lambda(x)\}, \qquad w_\lambda=(u-u_\lambda)^+\mathbf1_{\Sigma_\lambda}.
\end{aligned}
    \]
    Testing \(-\Delta u+u=(I_\alpha*u^p)u^{p-1}\) against \(\phi^2/u\) gives
\begin{equation}
  \label{eq:intro-wei-sob}
   \int_{\mathbb R^N}(I_\alpha*u^p)u^{p-2}\phi^2\,dx \le \|\phi\|_{H^1(\mathbb R^N)}^2,
   \qquad \phi\in H^1(\mathbb R^N);
 \end{equation}
 see also~\cite[Lemma~2.4]{PT06}. \eqref{eq:intro-wei-sob} holds with
  \(u_\lambda\) in place of \(u\). Testing the difference of the equations for \(u\)
  and \(u_\lambda\) with \(w_\lambda\) yields
  \[
  \begin{aligned}
 (2-p)\|w_\lambda\|_{H^1(\R^N)}^2 &\le p\int_{\mathbb R^N} \bigl(I_\alpha*(u^{p-1}w_\lambda)\bigr) u^{p-1}w_\lambda\,dx \le
  C\|u\|_{L^{2Np/(N+\alpha)}(E_\lambda)}^{2p-2} \|w_\lambda\|_{H^1(\R^N)}^2.
   \end{aligned}
   \]
 Thus \(w_\lambda=0\) whenever the coefficient on the right is smaller
 than \(2-p\). Integrability of \(u\) gives this condition for sufficiently negative
 \(\lambda\). At a stopping plane \(\lambda_0\), strict inequality
 \(u<u_{\lambda_0}\) in \(\Sigma_{\lambda_0}\) would imply, by dominated
 convergence, that the same coefficient tends to zero as \(\lambda\downarrow\lambda_0\).
 The plane could then move further. Hence \(u=u_{\lambda_0}\). Applying the argument in
 every direction, and using Hopf's lemma, gives radial symmetry and
 strict radial decrease.

  \medskip
 \noindent \emph{Morse index and nondegeneracy.} For an arbitrary positive
 solution, Morse index one cannot be inferred from a minimizing characterization. We
 prove that there exists \(c_Q>0\) such that
 \[
    \langle h,L_Qh\rangle \ge c_Q\|h\|_{H^1(\mathbb R^N)}^2, \quad h\in H^1_{\mathrm{rad}}(\mathbb R^N) \quad
    \text{satisfying} \int_{\mathbb R^N} (I_\alpha*Q^p)Q^{p-1}h\,dx=0.
    \]
   After centering \(Q\), let \(A_0=-\Delta+1-(I_\alpha*Q^p)Q^{p-2}\) and
   \(d\nu=(I_\alpha*Q^p)Q^p\,dx\). We work on \(T=\{Q\eta\in
   H^1_{\mathrm{rad}}(\mathbb R^N): \int_{\mathbb R^N}\eta\,d\nu=0\}\),
   where \(A_0\) is coercive.

We prove in Lemma~\ref{lem:derivative-modulus} that, for every \(Q\eta\in T\),
there exists \(Q\widetilde\eta\in T\) with \(\widetilde\eta'=|\eta'|\) such that
   \begin{equation}
  \label{eq:intro-compare}
    \langle Q\widetilde\eta,A_0(Q\widetilde\eta)\rangle =\langle Q\eta,A_0(Q\eta)\rangle,\qquad
    \langle Q\widetilde\eta,(A_0-L_Q)(Q\widetilde\eta)\rangle \ge\langle
    Q\eta,(A_0-L_Q)(Q\eta)\rangle.
   \end{equation}
  Since \(A_0-L_Q:H^1_{\mathrm{rad}}(\R^N)\to H^{-1}_{\mathrm{rad}}(\R^N)\) is compact, the supremum
    \[
 \tau=\sup_{h\in T\setminus\{0\}} \frac{\langle h,(A_0-L_Q)h\rangle} {\langle h,A_0h\rangle}>0
 \]
is attained. \eqref{eq:intro-compare} gives a maximizer \(h=Q\eta\) with
  \(\eta'\ge0\). Its Euler--Lagrange equation is \(\tau A_0h=(A_0-L_Q)h\).
Combining this with \(L_Q\left(\frac{\alpha+2}{4(p-1)}Q+\frac12x\cdot\nabla
  Q\right)=-Q\) gives
\[
   \underbrace{\int_{\mathbb R^N}Q^2\eta\,dx}_{>0} = \frac12\left(1-\frac1\tau\right) \underbrace{ \int_{\mathbb
     R^N}\eta Q^2x\cdot\nabla \bigl((I_\alpha*Q^p)Q^{p-2}\bigr)\,dx }_{<0\quad (\text{See
     Lemma~}\ref{lem:radial-virial-signs})}.
 \]
    Consequently, \(0<\tau<1\), and
\[
  \langle h,L_Qh\rangle \ge(1-\tau)\langle h,A_0h\rangle \ge c_Q\|h\|_{H^1(\R^N)}^2,
  \qquad h\in T.
    \]
 The identity \(L_QQ=-2(p-1)(I_\alpha*Q^p)Q^{p-1}\) provides a negative direction and places every
 radial kernel element in \(T\). Hence the radial Morse index is one and the radial kernel is
 trivial.

 The other sectors are treated by spherical harmonics, as in~\cite{Lenzmann09,Chen21}. The angular Riesz kernels
 obtained from the heat-semigroup representation~\cite{Stinga19,BM24} are positive and decrease with the
 harmonic degree. Their comparison shows that the first sector is nonnegative with kernel generated
 by \(Q'\), while all higher sectors are positive. For \(N=1\), the corresponding comparison applies to the odd
 subspace. This gives \(\operatorname{ind}L_Q=1\) and \(\ker L_Q=\operatorname{span}
 \{\partial_1Q,\ldots,\partial_NQ\}\).

 \medskip
  \noindent \emph{Uniqueness.} For fixed \(2<p<2^*\), we begin with the limit
  \(\alpha\uparrow N\), studied for ground states in~\cite[Sections~4--5]{Seok19}.
   Set \(\gamma=N-\alpha\) and \(U=c_{N,\alpha}^{1/(2p-2)}Q\). Then
\[
   -\Delta U+U=(|\cdot|^{-\gamma}*U^p)U^{p-1}.
  \]
    We establish \(c\le U(0)\le C\) and \(U(x)+|\nabla U(x)|\le Ce^{-c|x|}\) uniformly for all centered positive
      solutions with sufficiently small \(\gamma\). The bounds follow from a blow-up argument using \cite{GS81}.
      They imply convergence in \(H^1(\R^N)\cap L^\infty(\R^N)\) to
\[
U_*= \left(\int_{\mathbb R^N}w^p\,dx\right)^{-1/(2p-2)}w, \qquad -\Delta w+w=w^{p-1}, \qquad w>0,\quad w\in
H^1_{\mathrm{rad}}(\mathbb R^N),
\]
    where uniqueness of \(w\) follows from~\cite{Kwong89}. Thus the same limit is
    obtained for every sequence of centered positive solutions.

  At this limit, the linearized operator is
    \[
L_*h= \bigl(-\Delta+1-(p-1)w^{p-2}\bigr)h - \frac{p\,w^{p-1}}{\displaystyle\int_{\mathbb R^N}w^p\,dx}
\int_{\mathbb R^N}w^{p-1}h\,dx.
    \]
  The local operator is radially invertible by~\cite[Lemma~2.1]{CGNT07} and maps \(w\) to
   \(-(p-2)w^{p-1}\). These two facts give radial invertibility of \(L_*\). Convergence of
   the linearized operators then excludes distinct solutions approaching \(U_*\) and yields uniqueness
  for \(\alpha\) sufficiently close to \(N\).

  To continue in \(\alpha\), we use Li's a priori estimates~\cite[Theorem~1.5]{Li22}.
    They give compactness of centered positive solutions over compact subintervals of
    \(J_p=(\max\{0,(N-2)p-N\},N)\). Radial nondegeneracy provides a locally
  unique positive branch through each solution. Compactness therefore makes the number
    of centered positive solutions finite and locally constant on \(J_p\);
  compare~\cite[Section~8]{FLS16}. Since this number equals one near \(N\), it equals one
    throughout \(J_p\).

Finally, at \(p=2\), radial nondegeneracy allows us to continue each positive
 profile to exponents \(p>2\) by a parameter-dependent implicit function theorem in
 \(H^1_{\mathrm{rad}}(\R^N)\cap L^\infty(\R^N)\).

  The paper is organized as follows. Section~\ref{sec:symmetry-section} establishes radial symmetry.
    Section~\ref{sec:tech} provides the decay and quadratic-form estimates used in Section~\ref{sec:nondegeneracy}
    to prove Morse index one and nondegeneracy. In Section~\ref{sec:uniqueness}, we prove uniqueness.

\subsection*{Notation}

 Unless otherwise stated, all function spaces are over the reals.

\begin{itemize}
\item
  We use the usual notation for Lebesgue and Sobolev spaces.

\item
    We set \(H^{-1}(\R^N)=(H^1(\R^N))^*\). The brackets \(\langle\cdot,\cdot\rangle\) denote either
    the real \(L^2(\R^N)\) inner product or the duality pairing, with the convention
    \[
\langle f,g\rangle=\int_{\R^N}fg\,dx, \qquad \langle h,F\rangle=F[h] \quad
 \bigl(h\in H^1(\R^N),\ F\in H^{-1}(\R^N)\bigr).
  \]
   For \(h\in H^1(\R^N)\), the expression \(\langle h,L_Qh\rangle\) denotes the quadratic form of \(L_Q\).
   Operator inequalities are understood in the sense of quadratic forms.

\item
    We write
    \[
   B_R=\{x\in\R^N:|x|<R\},\qquad \mathbb S^{N-1}=\{x\in\R^N:|x|=1\},
  \]
and denote the surface measure of the unit sphere by \(|\mathbb S^{N-1}|\), with \(|\mathbb S^0|=2\). The
  inner product on the radial profile space is
\[
   \langle f,g\rangle_{L^2((0,\infty),r^{N-1}\,dr)} =\int_0^\infty f(r)g(r)r^{N-1}\,dr.
  \]

\item
  For a self-adjoint operator \(L\), the symbols \(\sigma(L)\) and \(\sigma_{\mathrm{ess}}(L)\)
  denote its spectrum and essential spectrum. The Morse index of \(L_Q\) is
 \[
\operatorname{ind}L_Q =\sup\left\{ \dim E:
   \begin{array}{l}
    E\subset H^1(\R^N)\text{ is a linear subspace},\\
 \langle h,L_Qh\rangle<0 \text{ for every }h\in E\setminus\{0\}
  \end{array}
   \right\}.
  \]
 Equivalently, it is the number of negative eigenvalues, counted with multiplicity. The radial Morse
   index \(\operatorname{ind}_{\mathrm{rad}}L_Q\) is defined by restricting \(E\) to
   subspaces of \(H^1_{\mathrm{rad}}(\R^N)\). The symbol \(\oplus\) denotes a direct sum.

\item
 For a measurable set \(E\), \(\mathbf1_E\) denotes its indicator function. For a real-valued function \(f\), we
 set
    \[
  f_+=\max\{f,0\},\qquad f_-=\max\{-f,0\}.
\]

\item
    The letters \(C\) and \(c\) denote positive constants whose values may change from one occurrence to the next.
     Dependence on varying parameters is stated where relevant. The notation \(o(1)\) denotes a quantity tending to
     zero in the limit under consideration; any required uniformity is specified explicitly.
\end{itemize}

\section{Symmetry}\label{sec:symmetry-section}
 The following regularity result allows us to compare a positive
  solution with its reflections throughout \eqref{eq:existence-range}.

\begin{lemma}[{\cite[Proposition~4.1]{MV13}}]\label{lem:regularity}
Let \( N\ge1\), \(0<\alpha<N\),  \(p\ge2\), \((N-2)p<N+\alpha\). Then every positive solution \(Q\in H^1(\R^N)\) of
  \eqref{eq:equation} satisfies
\[
 Q\in L^1(\R^N)\cap C^\infty(\R^N)\cap W^{2,r}(\R^N) \qquad (1<r<\infty).
  \]
\end{lemma}

\subsection{Proof of Theorem~\ref{thm:radial-symmetry}}

\begin{proof}[Proof of Theorem~\ref{thm:radial-symmetry}]
 For \(p\ge2\), this follow from \cite[Theorem~1.1]{BSR26} to the pair \((u,u)\), with equal
   exponents and \(\tau=\eta=1\).

It remains to consider \( 1+\frac{\alpha}{N}<p<2\). We use the method of moving planes~\cite{GNN79,CLO06}.

\smallskip
  \noindent\emph{Step 1.} We first derive the comparison estimates. By
   lemma~\ref{lem:regularity}, we can known
   \[
 u\in L^\infty(\R^N), \qquad \lim_{|x|\to\infty}u(x)=0, \qquad I_\alpha*u^p\in L^\infty(\R^N)\cap C(\R^N).
  \]
 Fix \(e\in\mathbb S^{N-1}\). For \(\lambda\in\R\), set
\[
   \Sigma_\lambda=\{x\in\R^N:x\cdot e<\lambda\}, \qquad x^\lambda=x+2(\lambda-x\cdot e)e,
  \]
    \[
   u_\lambda(x)=u(x^\lambda), \qquad E_\lambda=\{x\in\Sigma_\lambda:u(x)>u_\lambda(x)\}.
\]
    Since \(u-u_\lambda\) has zero trace on \(\partial\Sigma_\lambda\),
    \(w_\lambda=(u-u_\lambda)^+\mathbf1_{\Sigma_\lambda} \in H^1(\R^N)\). Then, we
    obtain
 \begin{equation}
\label{eq:reflection}
  \begin{aligned}
   (I_\alpha*u^p)(x^\lambda) &=\int_{\R^N}I_\alpha(x^\lambda-z^\lambda) u(z^\lambda)^p\,dz
   =\int_{\R^N}I_\alpha(x-z)u_\lambda(z)^p\,dz =(I_\alpha*u_\lambda^p)(x).
    \end{aligned}
  \end{equation}
Moreover, for \(x,y\in\Sigma_\lambda\),
\begin{equation}
  \label{eq:sigma-lmabda-positive}
  |x^\lambda-y|^2-|x-y|^2 =4(\lambda-x\cdot e)(\lambda-y\cdot e)>0.
\end{equation}
 By \eqref{eq:reflection}, \eqref{eq:sigma-lmabda-positive},
\begin{equation}\label{eq:symmetry-reflection}
    \begin{aligned}
(I_\alpha*u^p)(x)-(I_\alpha*u_\lambda^p)(x) &= \int_{\Sigma_\lambda}
 \bigl[I_\alpha(x-y)-I_\alpha(x^\lambda-y)\bigr] \bigl[u(y)^p-u_\lambda(y)^p\bigr]\,dy\\
 &\le \int_{E_\lambda} \bigl[I_\alpha(x-y)-I_\alpha(x^\lambda-y)\bigr] \bigl[u(y)^p-u_\lambda(y)^p\bigr]\,dy\\
   &\le p\bigl(I_\alpha*(u^{p-1}w_\lambda)\bigr)(x), \qquad x\in\Sigma_\lambda.
\end{aligned}
  \end{equation}

By the HLS, H\"older and Sobolev inequalities,
   \begin{equation}\label{eq:symmetry-hls}
\begin{aligned}
    \int_{\R^N} \bigl(I_\alpha*(u^{p-1}w_\lambda)\bigr) u^{p-1}w_\lambda\,dx &\le
    C\|u^{p-1}w_\lambda\|_{L^{2N/(N+\alpha)}(\R^N)}^2\\
 &\le C\|u\|_{L^{2Np/(N+\alpha)}(E_\lambda)}^{2p-2} \|w_\lambda\|_{H^1(\R^N)}^2,
    \end{aligned}
  \end{equation}
    where \(C>0\) depends only on \(N\), \(\alpha\), and \(p\).

    Let \(\phi\in C_c^\infty(\R^N)\) be real. Since \(u>0\) and smooth,
     \(\phi^2/u\in H^1(\R^N)\) and has compact support. Testing
  \begin{equation}
\label{eq:u}
  -\Delta u+u=(I_\alpha*u^p)u^{p-1}
\end{equation}
  with \(\phi^2/u\), equivalently by applying \cite[Lemma~2.4]{PT06}, we
  obtain
  \begin{equation*}
 \int_{\R^N}\left( |\nabla\phi|^2+\phi^2 -(I_\alpha*u^p)u^{p-2}\phi^2 \right)\,dx =
   \int_{\R^N} u^2\left|\nabla\left(\frac{\phi}{u}\right)\right|^2\,dx \ge0.
   \end{equation*}
 Approximating in \(H^1(\R^N)\) and almost everywhere by functions in \(C_c^\infty(\R^N)\), and applying
  Fatou's lemma, we deduce that
\begin{equation}\label{eq:symmetry-form-bound}
   \int_{\R^N}(I_\alpha*u^p)u^{p-2}\phi^2\,dx \le \|\phi\|_{H^1(\R^N)}^2 \qquad \bigl(\phi\in
    H^1(\R^N)\bigr).
  \end{equation}
  The same estimate holds with \(u\) replaced by \(u_\lambda\).

    Since both \(u\) and \(u_\lambda\) solve \eqref{eq:u},
    subtracting and taking \(w_\lambda\) as a test function, together
     with \eqref{eq:symmetry-reflection}, \eqref{eq:symmetry-hls}
     and \eqref{eq:symmetry-form-bound} gives
    \[
   \begin{aligned}
  \|w_\lambda\|_{H^1(\R^N)}^2 &= \int_{E_\lambda} \left[ (I_\alpha*u^p)u^{p-1} -(I_\alpha*u_\lambda^p)u_\lambda^{p-1}
  \right]w_\lambda\,dx\\
    &=\int_{E_\lambda}(I_\alpha*u_\lambda^p) \bigl(u^{p-1}-u_\lambda^{p-1}\bigr)w_\lambda\,dx
     +\int_{E_\lambda}\bigl[(I_\alpha*u^p)-(I_\alpha*u_\lambda^p)\bigr] u^{p-1}w_\lambda\,dx\\
  &\le (p-1)\int_{E_\lambda} (I_\alpha*u_\lambda^p) u_\lambda^{p-2}w_\lambda^2\,dx+
  p\int_{E_\lambda} \bigl(I_\alpha*(u^{p-1}w_\lambda)\bigr) u^{p-1}w_\lambda\,dx\\
    &\le (p-1)\|w_\lambda\|_{H^1(\R^N)}^2+ C\|u\|_{L^{2Np/(N+\alpha)}(E_\lambda)}^{2p-2}
     \|w_\lambda\|_{H^1(\R^N)}^2.
    \end{aligned}
\]
   Thus,
    \begin{equation}\label{eq:symmetry-subquadratic}
 (2-p)\|w_\lambda\|_{H^1(\R^N)}^2 \le C\|u\|_{L^{2Np/(N+\alpha)}(E_\lambda)}^{2p-2} \|w_\lambda\|_{H^1(\R^N)}^2.
    \end{equation}

\smallskip
\noindent\emph{Step 2.} We show that \(u\le u_\lambda\) in
 \(\Sigma_\lambda\) for all sufficiently negative \(\lambda\). Since
 \(E_\lambda\subset\Sigma_\lambda\) and
 \(\|u\|_{L^{2Np/(N+\alpha)}(\Sigma_\lambda)}\to0\) as
 \(\lambda\to-\infty\), estimate~\eqref{eq:symmetry-subquadratic} gives
 \(w_\lambda=0\) for all sufficiently negative \(\lambda\). By
 continuity, \(u\le u_\lambda\) in \(\Sigma_\lambda\) for all sufficiently
 negative \(\lambda\).

\smallskip
 \noindent\emph{Step 3.} We move the plane. Set
\[
  \lambda_0= \sup\left\{ \lambda\in\R: u\le u_\mu\text{ in }\Sigma_\mu \text{ for every }\mu\le\lambda \right\}.
\]
Step 2 shows that \(\lambda_0>-\infty\). If \(\lambda_0=\infty\), then, for
each fixed \(x\in\R^N\), \(0<u(x)\le u(x^\lambda)\to0\) as
\(\lambda\to\infty\), a contradiction. Thus \(\lambda_0\in\R\). By continuity and
\eqref{eq:symmetry-reflection},
  \[
  u\le u_{\lambda_0}, \qquad I_\alpha*u^p\le I_\alpha*u_{\lambda_0}^p \quad\text{in
   }\Sigma_{\lambda_0}.
    \]
  Consequently,
    \[
 (-\Delta+1)(u_{\lambda_0}-u) =(I_\alpha*u_{\lambda_0}^p)u_{\lambda_0}^{p-1} -(I_\alpha*u^p)u^{p-1}
   \ge0 \quad\text{in }\Sigma_{\lambda_0}.
    \]
 The strong maximum principle implies that either \(u_{\lambda_0}\equiv u\) or \(u_{\lambda_0}>u\) in
 \(\Sigma_{\lambda_0}\).

   Suppose that \(u_{\lambda_0}>u\) in \(\Sigma_{\lambda_0}\). For every fixed \(x\)
   with \(x\cdot e<\lambda_0\), continuity gives \(u(x)<u_\lambda(x)\) when
   \(\lambda>\lambda_0\) is sufficiently close to \(\lambda_0\). For
   \(x\cdot e>\lambda_0\), the point \(x\) lies outside
   \(\Sigma_\lambda\) for all such \(\lambda\). Thus
   \(\mathbf1_{E_\lambda}\to0\) almost everywhere as
   \(\lambda\downarrow\lambda_0\). By dominated convergence,
\[
  \|u\|_{L^{2Np/(N+\alpha)}(E_\lambda)}\to0 \qquad(\lambda\downarrow\lambda_0).
  \]
By~\eqref{eq:symmetry-subquadratic}, \(w_\lambda=0\) for every
\(\lambda\in(\lambda_0,\lambda_0+\varepsilon]\), for some \(\varepsilon>0\), contradicting the
 definition of \(\lambda_0\). Hence \(u=u_{\lambda_0}\) in \(\R^N\).

\smallskip
\noindent\emph{Step 4.} It remains to identify the center and to obtain strict radial decrease.
For every \(\lambda<\lambda_0\), the strong maximum principle gives \(u_\lambda>u\) in
\(\Sigma_\lambda\), unless \(u=u_\lambda\). The latter alternative is impossible: symmetry about
the two parallel planes would imply
 \[
u(x)=u\bigl(x+2k(\lambda_0-\lambda)e\bigr) \to0 \qquad(k\to\infty),
    \]
  contrary to \(u>0\). Hopf's lemma on \(\partial\Sigma_\lambda\) gives
    \(-2\partial_eu=\partial_e(u_\lambda-u)<0\). Together with reflection about \(x\cdot e=\lambda_0\),
    this yields
   \[
  \partial_eu(x)>0\quad(x\cdot e<\lambda_0),\qquad \partial_eu(x)<0\quad(x\cdot e>\lambda_0).
    \]
    The direction \(e\) was arbitrary. Since \(u>0\) and \(u(x)\to0\), it attains its maximum at some
    \(x_0\in\R^N\). The equality \(\nabla u(x_0)=0\) forces \(\lambda_0(e)=x_0\cdot e\) for every \(e\in\mathbb
    S^{N-1}\). Thus \(u\) is invariant under every reflection through \(x_0\), so \(u(x)=U(|x-x_0|)\).
    Finally, \(U'(r)=-\partial_eu(x_0-re)<0\) for \(r>0\).
\end{proof}

\section{Technical estimates}\label{sec:tech}

 We now assume~\eqref{eq:range}, so that \(p\ge2\). We record the decay and spectral
 properties needed below, centering each positive solution at the origin by
 Theorem~\ref{thm:radial-symmetry}.

 \begin{lemma}\label{lem:nonlocal-quadratic-bound}
Let \( N\ge1\), \(0<\alpha<N\),  \(p\ge2\), \((N-2)p<N+\alpha\), and let \(Q\) be a positive solution of \eqref{eq:equation}. Then for
  \(h\in L^2(\R^N)\),
 \begin{equation*}
    \left| \int_{\R^N} (I_\alpha*(Q^{p-1}h))Q^{p-1}h\,dx \right| \le \int_{\R^N}
    (I_\alpha*Q^p)Q^{p-2}h^2\,dx.
   \end{equation*}
 \end{lemma}

\begin{proof}
  Since \(I_\alpha>0\),
  \[
  \left| \int_{\R^N} (I_\alpha*(Q^{p-1}h))Q^{p-1}h\,dx \right|\le \iint_{\R^N\times\R^N} I_\alpha(x-y)
   Q(x)^{p-1}Q(y)^{p-1} |h(x)h(y)|\,dx\,dy.
\]
 Using
  \[
Q(x)^{p-1}Q(y)^{p-1}|h(x)h(y)|\le \frac12 Q(x)^{p-2}h(x)^2Q(y)^p + \frac12 Q(y)^{p-2}h(y)^2Q(x)^p,
    \]
 we obtain, by the symmetry \(I_\alpha(x-y)=I_\alpha(y-x)\),
\[
  \left| \int_{\R^N} (I_\alpha*(Q^{p-1}h))Q^{p-1}h\,dx \right|\le \int_{\R^N} (I_\alpha*Q^p)Q^{p-2}h^2\,dx.
   \]
 This completes the proof.
\end{proof}

\begin{lemma}\label{lem:decay}
Let \( N\ge1\), \(0<\alpha<N\),  \(p\ge2\), \((N-2)p<N+\alpha\), and let \(Q\) be a positive solution of
    \eqref{eq:equation}. Then, after centering \(Q\) at the origin,
\[
 Q(x)+|\nabla Q(x)|\le Ce^{-|x|/2},
  \]
   and
 \[
(I_\alpha*Q^p)Q^{p-2}\in L^\infty(\R^N), \qquad \lim_{|x|\to\infty} (I_\alpha*Q^p)(x)Q(x)^{p-2}=0.
\]
Moreover, the operator
 \[
 h\longmapsto Q^{p-1}\bigl(I_\alpha*(Q^{p-1}h)\bigr)
\]
  is compact on \(L^2(\R^N)\), \(L_Q\) is self-adjoint on \(L^2(\R^N)\) with
   domain \(H^2(\R^N)\), and
   \(\sigma_{\mathrm{ess}}(L_Q)=[1,\infty)\).
\end{lemma}

   \begin{proof}
By Theorem~\ref{thm:radial-symmetry}, after a translation, \(Q\) is radial
and decreasing. By~\cite[Theorem~1.5]{Li22},
    \[
  Q(x)+|\nabla Q(x)|\le Ce^{-|x|/2}.
 \]
 By~\cite[Lemma~6.2]{MV13},
    \[
   (I_\alpha*Q^p)(x)=O\bigl(|x|^{\alpha-N}\bigr) \qquad (|x|\to\infty).
  \]
Since \(p\ge2\), it follows that
 \[
 (I_\alpha*Q^p)Q^{p-2}\in L^\infty(\R^N), \qquad (I_\alpha*Q^p)(x)Q(x)^{p-2}\to0 \quad (|x|\to\infty).
\]

 We next show that
\begin{equation}
   \label{eq:cpt-Qp-1}
h\longmapsto Q^{p-1}\bigl(I_\alpha*(Q^{p-1}h)\bigr)\quad \text{is compact on} \quad L^2(\R^N).
   \end{equation}
    By Cwikel's theorem~\cite{Cwikel77}, applied with exponent \(2N/\alpha>2\), the
    operator \(h\mapsto I_{\alpha/2}*(Q^{p-1}h)\) is compact on \(L^2(\R^N)\). Its adjoint
    \(h\mapsto Q^{p-1}(I_{\alpha/2}*h)\) is bounded on \(L^2(\R^N)\); indeed, Hölder's and
    HLS inequalities give
\[
\begin{aligned}
  \|Q^{p-1}(I_{\alpha/2}*h)\|_{L^2(\R^N)} &\le \|Q^{p-1}\|_{L^{2N/\alpha}(\R^N)} \|I_{\alpha/2}*h\|_{L^{2N/(N-\alpha)}(\R^N)}\\
    &\le C\|Q^{p-1}\|_{L^{2N/\alpha}(\R^N)}\|h\|_{L^2(\R^N)}.
    \end{aligned}
  \]
Since \(I_{\alpha/2}*I_{\alpha/2}=I_\alpha\), \eqref{eq:cpt-Qp-1} holds.

  Since \((I_\alpha*Q^p)Q^{p-2}\) is bounded and vanishes at infinity, \cite[Lemma~7.11]{Teschl09}
  gives
  \begin{equation}
\label{eq:cpt-Qp}
   \bigl[(I_\alpha*Q^p)Q^{p-2}\bigr](-\Delta+1)^{-1} \colon L^2(\R^N)\to L^2(\R^N) \quad\text{is compact}.
    \end{equation}

 Multiplication by \((I_\alpha*Q^p)Q^{p-2}\) is bounded and symmetric, as
  is the compact operator
  \(Q^{p-1}\bigl(I_\alpha*(Q^{p-1}\,\cdot)\bigr)\). Hence, by the
 Kato--Rellich theorem \cite[Theorem~6.4]{Teschl09}, \(L_Q\)
  is self-adjoint on \(L^2(\R^N)\) with domain \(H^2(\R^N)\).

 Moreover, by \eqref{eq:linearization}, \eqref{eq:cpt-Qp-1} and \eqref{eq:cpt-Qp},
  \(\bigl(L_Q-(-\Delta+1)\bigr)(-\Delta+1)^{-1}\) is compact on \(L^2(\R^N)\).
 Thus \(L_Q-(-\Delta+1)\) is relatively compact with respect to \(-\Delta+1\).
  By Weyl's theorem \cite[Theorem~6.19]{Teschl09},
  \[
   \sigma_{\mathrm{ess}}(L_Q) = \sigma_{\mathrm{ess}}(-\Delta+1) = [1,\infty).
  \]
 This completes the proof.
\end{proof}

\subsection{Riesz interactions of radial functions}
The comparison of the quadratic forms at \(Q\eta\) and \(Q\widetilde\eta\), where
\(\widetilde\eta'=|\eta'|\), uses the strict inequality \eqref{pge2:positive-covariance}. We
first compute \(\int_{B_R}(I_\alpha*\mathbf1_{B_S})\,dx\) in
Lemma~\ref{lem:ball-kernel}. Lemma~\ref{lem:association} then establishes
\eqref{pge2:positive-covariance}, which is applied to \(L_Q\) in
Lemma~\ref{lem:derivative-modulus}.

\begin{lemma}\label{lem:ball-kernel}
 For \(N\ge1\), \(0<\alpha<N\), and \(R,S>0\),
    \begin{equation}\label{pge2:ball-gram}
  \begin{aligned}
    \int_{B_R}(I_\alpha*\mathbf1_{B_S})(x)\,dx =C_{N,\alpha}\int_0^{\min(R^2,S^2)} t^{(N-\alpha)/2-1}
    (R^2-t)^{\alpha/2}(S^2-t)^{\alpha/2}\,dt,
\end{aligned}
\end{equation}
   where \( C_{N,\alpha} = \frac{2^{-\alpha}\pi^{N/2}} {\Gamma(N/2)\Gamma(1+\alpha/2)^2}. \)
\end{lemma}

\begin{proof}
  By~\cite[Theorem~1.2]{GCO23}, we obtain
  \[
(I_\alpha*\mathbf1_{B_S})(r) = \frac{c_{N,\alpha}|\mathbb S^{N-1}|}{\alpha}S^\alpha {}_2F_1\left(
  \frac{N-\alpha}{2},-\frac{\alpha}{2}; \frac N2;\frac{r^2}{S^2} \right), \qquad 0\le r<S.
\]
 Here \({}_2F_1\) is the Gauss hypergeometric function.

Assume first that \(0<R<S\). Using polar coordinates, \cite[Eqns.~(15.2.1)
and~(15.6.1)]{OLBC10}, and the change of variables \(t=R^2u\), we obtain
   \[
    \begin{aligned}
  \int_{B_R}(I_\alpha*\mathbf1_{B_S})(x)\,dx &= \frac{c_{N,\alpha}|\mathbb S^{N-1}|^2}{\alpha}S^\alpha \int_0^R r^{N-1}
   {}_2F_1\left( \frac{N-\alpha}{2},-\frac{\alpha}{2}; \frac N2;\frac{r^2}{S^2} \right)\,dr\\
  &= \frac{c_{N,\alpha}|\mathbb S^{N-1}|^2}{N\alpha} R^NS^\alpha {}_2F_1\left( \frac{N-\alpha}{2},-\frac{\alpha}{2}; \frac
    N2+1;\frac{R^2}{S^2} \right)\\
 &= \frac{c_{N,\alpha}|\mathbb S^{N-1}|^2}{N\alpha} \frac{\Gamma(N/2+1)} {\Gamma((N-\alpha)/2)\Gamma(1+\alpha/2)} R^NS^\alpha\\
&\qquad\cdot \int_0^1 u^{(N-\alpha)/2-1} (1-u)^{\alpha/2} \left(1-\frac{R^2}{S^2}u\right)^{\alpha/2}\,du\\
&= \frac{2^{-\alpha}\pi^{N/2}} {\Gamma(N/2)\Gamma(1+\alpha/2)^2} \int_0^{R^2}
t^{(N-\alpha)/2-1} (R^2-t)^{\alpha/2}(S^2-t)^{\alpha/2}\,dt.
 \end{aligned}
\]
    The case \(R=S\) follows by monotone convergence as \(R\uparrow S\). Symmetry in \(R\) and
    \(S\) then gives \eqref{pge2:ball-gram} for all \(R,S>0\).
\end{proof}

\begin{lemma}\label{lem:association}
  Let \(Q\) be a positive solution of~\eqref{eq:equation}, centered at the origin, and set
  \begin{equation}\label{eq:centered-measure}
  d\nu(x)=Q(x)^p(I_\alpha*Q^p)(x)\,dx,\qquad e_t(x)=\mathbf1_{\{|x|>t\}}-
  \frac{\nu(\{|x|>t\})}{\nu(\R^N)}.
   \end{equation}
   Then
    \begin{equation}\label{pge2:positive-covariance}
 \int_{\R^N} \bigl(I_\alpha*(Q^pe_s)\bigr)Q^pe_t\,dx>0 \qquad (s,t>0).
\end{equation}
\end{lemma}

\begin{proof}
 Since \(Q\) is radial and \(Q(r)\to0\) as \(r\to\infty\), the fundamental theorem of calculus
   gives, for \(r=|x|\),
 \begin{equation}
\label{eq:layer}
  Q(x)^p\mathbf1_{B_s}(x) = -\mathbf1_{\{r<s\}}\int_r^\infty (Q^p)'(R)\,dR =-\int_0^\infty
  \mathbf1_{B_{\min\{s,R\}}}(x)(Q^p)'(R)\,dR,
   \end{equation}
 where \(\min\{\infty,R\}=R\).

  By \eqref{pge2:ball-gram}, \eqref{eq:layer} and applying Tonelli's theorem, we obtain
\begin{equation}\label{eq:truncated-gram}
  \begin{aligned}
   &\int_{\R^N} \bigl(I_\alpha*(Q^p\mathbf1_{B_s})\bigr) Q^p\mathbf1_{B_t}\,dx\\
 &\quad= \int_0^\infty\int_0^\infty \left[ \int_{B_{\min\{t,T\}}} \bigl(I_\alpha*
   \mathbf1_{B_{\min\{s,R\}}}\bigr)(x)\,dx \right] (Q^p)'(R)(Q^p)'(T)\,dR\,dT\\
    &\quad=C_{N,\alpha} \int_0^\infty z^{(N-\alpha)/2-1} F_s(z)F_t(z)\,dz,
 \end{aligned}
    \end{equation}
where \[ F_s(z) = -\int_0^\infty (\min\{s,R\}^2-z)_+^{\alpha/2}(Q^p)'(R)\,dR, \qquad z>0.
  \]

   We next compare \(F_s\) with \(F_\infty\). For \(0<z<s^2\),
 \[
\begin{aligned}
   F_s(z) &= (s^2-z)^{\alpha/2} \left[ Q(s)^p -\int_0^s \left( \frac{(R^2-z)_+}{s^2-z} \right)^{\alpha/2} (Q^p)'(R)\,dR \right].
    \end{aligned}
 \]
 Moreover,
\[
\begin{aligned}
  F_\infty(z)-F_s(z) &= -(s^2-z)^{\alpha/2} \int_s^\infty \left[ \left( \frac{R^2-z}{s^2-z} \right)^{\alpha/2}-1
  \right](Q^p)'(R)\,dR.
  \end{aligned}
   \]
 Consequently,
   \[
   \frac{F_\infty(z)}{F_s(z)} = 1+ \frac{ -\displaystyle\int_s^\infty \left[
     \left(\frac{R^2-z}{s^2-z}\right)^{\alpha/2}-1 \right](Q^p)'(R)\,dR }{ Q(s)^p-\displaystyle\int_0^s \left(
     \frac{(R^2-z)_+}{s^2-z} \right)^{\alpha/2}(Q^p)'(R)\,dR }.
\]
  For \(R>s\) and \(0<z<s^2\), \(z\mapsto ((R^2-z)/(s^2-z))^{\alpha/2}\) is
   strictly increasing. On the other hand, if \(R<s\), then \(z\mapsto
   ((R^2-z)_+/(s^2-z))^{\alpha/2}\) is nonincreasing on \(0<z<s^2\): it is
   strictly decreasing for \(0<z<R^2\) and vanishes for \(R^2\le z<s^2\).
   Since \((Q^p)'(R)<0\), \(F_\infty(z)/F_s(z)\) is strictly increasing,
   and therefore \(F_s(z)/F_\infty(z)\) is strictly decreasing on
   \(0<z<s^2\), and it is zero for \(z\ge s^2\).

  By \(I_\alpha(x-y)=I_\alpha(y-x)\), Tonelli's theorem, taking \(t=\infty\)
  in~\eqref{eq:truncated-gram},
    \begin{equation}
 \label{eq:nuBs}
  \begin{aligned}
    \nu(B_s) &= \int_{\R^N} (I_\alpha*Q^p)Q^p\mathbf1_{B_s}\,dx= \int_{\R^N}
    \bigl(I_\alpha*(Q^p\mathbf1_{B_s})\bigr)Q^p\,dx\\
 &=C_{N,\alpha} \int_0^\infty z^{(N-\alpha)/2-1}F_s(z)F_\infty(z)\,dz.
\end{aligned}
 \end{equation}
 Taking \(s=t=\infty\) in~\eqref{eq:truncated-gram} gives
 \begin{equation}
    \label{eq:nurn}
 \nu(\R^N) = C_{N,\alpha} \int_0^\infty z^{(N-\alpha)/2-1}F_\infty(z)^2\,dz.
\end{equation}

  Using \eqref{eq:truncated-gram}, \eqref{eq:nuBs},
    \eqref{eq:nurn} and \(e_s = \frac{\nu(B_s)}{\nu(\R^N)}
    -\mathbf1_{B_s}\) a.e.,
\[
  \begin{aligned}
\int_{\R^N}(I_\alpha*(Q^pe_s))Q^pe_t\,dx &= C_{N,\alpha}\int_0^\infty
z^{(N-\alpha)/2-1}F_s(z)F_t(z)\,dz -\frac{\nu(B_s)\nu(B_t)}{\nu(\R^N)}\\
  &= \frac{C_{N,\alpha}^2}{2\nu(\R^N)} \int_0^\infty\int_0^\infty (zw)^{(N-\alpha)/2-1} \bigl[ F_s(z)F_\infty(w)-F_s(w)F_\infty(z) \bigr]\\
  &\qquad\qquad\cdot \bigl[ F_t(z)F_\infty(w)-F_t(w)F_\infty(z) \bigr] \,dz\,dw.
  \end{aligned}
 \]
  If \(z<w\), the monotonicity of \(F_s/F_\infty\) and \(F_t/F_\infty\) gives
 \[
 F_s(z)F_\infty(w)-F_s(w)F_\infty(z)\ge0, \qquad F_t(z)F_\infty(w)-F_t(w)F_\infty(z)\ge0.
\]
 Both inequalities are strict whenever \(0<z<w<\min\{s^2,t^2\}\). For \(z>w\), both factors change sign,
   so their product remains nonnegative. Hence the integrand is nonnegative everywhere and is strictly positive on
   a set of positive measure. Therefore \eqref{pge2:positive-covariance} hold.
\end{proof}

\section{Nondegeneracy}\label{sec:nondegeneracy}

  Set
   \[
  A_0=-\Delta+1-(I_\alpha*Q^p)Q^{p-2}.
 \]
We use the fixed radial constraint space
 \[
T=\left\{h\in H^1_{\mathrm{rad}}(\R^N): \int_{\R^N}(I_\alpha*Q^p)Q^{p-1}h\,dx=0\right\}.
    \]
 The constraint defining \(T\) is determined by
 \[
   L_QQ=-2(p-1)(I_\alpha*Q^p)Q^{p-1}.
  \]
  Indeed, self-adjointness implies that every radial element of \(\ker L_Q\) belongs to \(T\). We shall prove in
  Proposition~\ref{prop:radial-coercivity} that the quadratic form of \(L_Q\) is coercive on \(T\). This will
  exclude radial kernel elements and, since \(\langle Q,L_QQ\rangle<0\) and \(T\) has codimension one
  in \(H^1_{\mathrm{rad}}(\R^N)\), establish Morse index one among radial perturbations.

   We begin with the quadratic form of \(A_0\). Its coercivity
   on \(T\), together with the compactness of the form of
   \(A_0-L_Q\), will allow us to attain the supremum
   in~\eqref{eq:constrained-maximum}.

\subsection{Radial comparison lemmas}

   For \(h=Q\eta\in C_c^\infty(\R^N)\), the ground-state representation \cite[Lemma~2.4]{PT06} gives
    \begin{equation}\label{eq:A0-ground-form}
 \langle Q\eta,A_0(Q\eta)\rangle =\int_{\R^N}Q(x)^2|\nabla\eta(x)|^2\,dx.
  \end{equation}
  Since \((I_\alpha*Q^p)Q^{p-2}\in L^\infty(\R^N)\), the quadratic form of \(A_0\) is continuous
    on \(H^1(\R^N)\). To extend \eqref{eq:A0-ground-form}, choose \(h_j\in C_c^\infty(\R^N)\)
    with \(h_j\to h\) in \(H^1(\R^N)\). Applying \eqref{eq:A0-ground-form} to \(h_j-h_k\) shows
    that \(Q\nabla(h_j/Q)\) is Cauchy in \(L^2(\R^N)\). Since \(Q\) is smooth and positive,
    its limit is \(Q\nabla(h/Q)\) on every compact set. Passing to the limit therefore proves
  \eqref{eq:A0-ground-form} for every \(h=Q\eta\in H^1(\R^N)\). Consequently,
\[
   A_0\ge0,\qquad \ker A_0=\operatorname{span}\{Q\}.
  \]
  \begin{lemma}
  \label{lem:ground-state-coercivity}
 There exists \(c_Q>0\) such that
 \begin{equation}
  \label{eq:coercivity-H1}
    \langle h,A_0h\rangle\ge c_Q\|h\|_{H^1(\R^N)}^2 \qquad(h\in T).
   \end{equation}
 Moreover, if \(h_j\rightharpoonup0\) in \(H^1(\R^N)\), then
  \begin{equation}
  \label{cpt-A0-Q0}
    \langle h_j,(A_0-L_Q)h_j\rangle\to0.
  \end{equation}
 \end{lemma}

\begin{proof}
 Equation~\eqref{eq:equation} gives \(A_0Q=0\), while \eqref{eq:linearization} gives
 \begin{equation}\label{eq:A0-minus-LQ}
A_0-L_Q=(p-2)(I_\alpha*Q^p)Q^{p-2} +pQ^{p-1}I_\alpha*(Q^{p-1}\,\cdot).
\end{equation}
 Since \((I_\alpha*Q^p)Q^{p-2}\) is bounded and tends to zero at infinity,
  its quadratic form is compact on \(H^1(\R^N)\) by local Rellich
  compactness. Suppose that \eqref{eq:coercivity-H1} fails. Then there
  exist \(h_j\in T\) such that
   \[
 \|h_j\|_{H^1(\R^N)}=1,\qquad \langle h_j,A_0h_j\rangle\to0,\qquad h_j\rightharpoonup h\quad\text{in
 }H^1(\R^N).
   \]
  The space \(T\) is weakly closed, and compactness gives
\[
   \int_{\R^N}(I_\alpha*Q^p)Q^{p-2}h^2\,dx =\lim_{j\to\infty} \int_{\R^N}(I_\alpha*Q^p)Q^{p-2}h_j^2\,dx=1.
 \]
   Weak lower semicontinuity then yields \(0\le\langle
   h,A_0h\rangle\le0\). Hence \(h\in\ker A_0\cap T=\{0\}\), a
   contradiction. This proves \eqref{eq:coercivity-H1}.

 By Lemma~\ref{lem:decay}, for \(h_j\rightharpoonup0\) in \(H^1(\R^N)\),
    \begin{equation}
      \label{eq:cpt-hj}
  \int_{\R^N}(I_\alpha*Q^p)Q^{p-2}h_j^2\,dx\to0, \qquad \int_{\R^N}(I_\alpha*(Q^{p-1}h_j))Q^{p-1}h_j\,dx \to0.
     \end{equation}
  The second limit  follows from the first by Lemma~\ref{lem:nonlocal-quadratic-bound}. Together with
   \eqref{eq:A0-minus-LQ}, \eqref{eq:cpt-hj} prove \eqref{cpt-A0-Q0}.
   \end{proof}

 \begin{lemma}
\label{lem:derivative-modulus}
    Let \(h=Q\eta\in T\). There exists a radial function \(\widetilde\eta\) such
    that
 \[
 Q\widetilde\eta\in T,\qquad \widetilde\eta'=|\eta'|\quad\text{a.e. on }(0,\infty),
  \]
  and
 \begin{equation}
    \label{eq:sign-increases}
   \langle Q\widetilde\eta,L_Q(Q\widetilde\eta) \rangle \le \langle Q\eta,L_Q(Q\eta) \rangle,\, \text{with strict inequality if } \,
    \eta'\, \text{ changes sign}.
  \end{equation}
\end{lemma}

\begin{proof}
 With \(d\nu\) and \(e_t\) as in~\eqref{eq:centered-measure}, \eqref{eq:A0-ground-form} and
 \eqref{eq:A0-minus-LQ} give
  \begin{equation}\label{eq:ground-form}
   \langle Q\eta,L_Q(Q\eta)\rangle =\int_{\R^N}Q^2|\nabla\eta|^2\,dx
   -(p-2)\int_{\R^N}\eta^2\,d\nu(x)-p\int_{\R^N}(I_\alpha*(Q^p\eta))Q^p\eta\,dx.
 \end{equation}
 By Lemma~\ref{lem:association}, for \(r,t>0\),
\begin{equation}\label{eq:positive-covariance}
 \int_{\R^N}e_re_t\,d\nu(x) =\frac{\nu(B_{\min(r,t)})\nu(\R^N\setminus B_{\max(r,t)})} {\nu(\R^N)}>0,\quad
   \int_{\R^N}(I_\alpha*(Q^pe_r))Q^pe_t\,dx>0.
    \end{equation}
    
 For \(R>2\), define
 \[
   \eta_R(r) =\int_1^r\mathbf1_{[R^{-1},R]}(s)\eta'(s)\,ds -\frac{1}{\nu(\R^N)} \int_{\R^N}\int_1^{|x|}
    \mathbf1_{[R^{-1},R]}(s)\eta'(s)\,ds\,d\nu(x),
    \]
   and define \(\widetilde\eta_R\) by replacing \(\eta'\) with \(|\eta'|\). These functions are
   bounded and have zero \(\nu\)-mean. So \(Q\eta_R,Q\widetilde\eta_R\in T\).

  By Lemma~\ref{lem:ground-state-coercivity} and \eqref{eq:A0-ground-form},
  \[
    \begin{aligned}
  c_Q\|Q(\eta_R-\eta)\|_{H^1(\R^N)}^2 &\le\int_{\R^N}Q^2 \bigl|\bigl(\mathbf1_{[R^{-1},R]}(|x|)-1\bigr)
   \eta'(|x|)\bigr|^2\,dx\to0,\\
   c_Q\|Q(\widetilde\eta_R-\widetilde\eta_S)\|_{H^1(\R^N)}^2 &\le\int_{\R^N}Q^2
     \bigl|\bigl(\mathbf1_{[R^{-1},R]}(|x|) -\mathbf1_{[S^{-1},S]}(|x|)\bigr)
     |\eta'(|x|)|\bigr|^2\,dx\to0
 \end{aligned}
\]
   as \(R,S\to\infty\). Thus \(Q\widetilde\eta_R\) converges in
   \(H^1(\R^N)\) to a function \(Q\widetilde\eta\in T\). Moreover,
   \begin{equation}\label{eq:def-title-eta}
 Q\widetilde\eta\in T,\qquad \widetilde\eta'=|\eta'|\quad\text{a.e. on }(0,\infty).
    \end{equation}

 If \(\eta'\) is compactly supported in \((0,\infty)\), Fubini's
   theorem and \(\int_{\R^N}\eta\,d\nu(x)=0\) give
    \[
\begin{aligned}
  \int_0^\infty e_t(x)\eta'(t)\,dt &= \int_0^{|x|}\eta'(t)\,dt - \frac{1}{\nu(\R^N)}
   \int_0^\infty\nu(\{|y|>t\})\eta'(t)\,dt\\
 &= \eta(|x|)-\eta(0) - \frac{1}{\nu(\R^N)} \int_{\R^N} \bigl(\eta(|y|)-\eta(0)\bigr)\,d\nu(y) = \eta(|x|).
    \end{aligned}
   \]
 Hence
    \begin{equation}\label{eq:covariance-form}
  \begin{aligned}
   &(p-2)\int_{\R^N}\eta^2\,d\nu(x) + p\int_{\R^N} \bigl(I_\alpha*(Q^p\eta)\bigr)Q^p\eta\,dx\\
   &\,= \iint_{(0,\infty)^2} \eta'(r)\eta'(t)\left[ (p-2)\int_{\R^N}e_re_t\,d\nu(x) + p\int_{\R^N}
   \bigl(I_\alpha*(Q^pe_r)\bigr)Q^pe_t\,dx \right]\,dr\,dt.
\end{aligned}
    \end{equation}
   Apply \eqref{eq:covariance-form} to \(\eta_R\) and \(\widetilde\eta_R\),
     whose derivatives are compactly supported. Using
     \(|\widetilde\eta_R'|=|\eta_R'|\) and \eqref{eq:ground-form}, we obtain
 \[
 \begin{aligned}
  &\langle Q\widetilde\eta_R,L_Q(Q\widetilde\eta_R)\rangle -\langle
  Q\eta_R,L_Q(Q\eta_R)\rangle\\
 &\quad= -4\iint_{[R^{-1},R]^2} (\eta')_+(r)(\eta')_-(t)\\
    &\qquad\qquad\cdot\left( (p-2)\int_{\R^N}e_re_td\nu(x) + p\int_{\R^N}(I_\alpha*(Q^pe_r))Q^pe_tdx
    \right)drdt.
\end{aligned}
\]
    The left-hand side converges by continuity of the quadratic form on
      \(H^1(\R^N)\). Monotone convergence on the right gives
   \[
\begin{aligned}
 &\langle Q\widetilde\eta,L_Q(Q\widetilde\eta)\rangle -\langle Q\eta,L_Q(Q\eta)\rangle\\
&\quad= -4\iint_{(0,\infty)^2} (\eta')_+(r)(\eta')_-(t)\\
&\qquad\qquad\cdot\left( (p-2)\int_{\R^N}e_re_td\nu(x) + p\int_{\R^N}(I_\alpha*(Q^pe_r))Q^pe_tdx \right)drdt.
  \end{aligned}
   \]
  This proves \eqref{eq:sign-increases}, including strictness.
\end{proof}

 Lemma~\ref{lem:derivative-modulus} allows a maximizer in \eqref{eq:constrained-maximum} to be chosen in
 the form \(Q\eta\), with \(\eta\) nondecreasing. The next lemma determines the signs of
  the two integrals that occur in~\eqref{eq:dilation-eigenvalue}. Together with the
  identity~\eqref{eq:frequency}, these signs will show that the supremum in
 \eqref{eq:constrained-maximum} is strictly smaller than one.

\begin{lemma}
 \label{lem:radial-virial-signs}
  Let \(Q\eta\in T\), where \(\eta\) is radial, nonconstant, and \(\eta'\ge0\) almost everywhere on
   \((0,\infty)\). Then
 \begin{equation}\label{eq:positive-unweighted-mean}
 \int_{\mathbb{R}^N}Q^2\eta\,dx>0.
\end{equation}
and
  \begin{equation}\label{eq:negative-dilation-weight}
  \int_{\R^N} \eta Q^2x\cdot\nabla\bigl((I_\alpha*Q^p)Q^{p-2}\bigr)\,dx<0.
\end{equation}
\end{lemma}

\begin{proof}
 For \(0<r<s\), let \(H=\{z\in\R^N:z_1<(r+s)/2\}\) and
 \(\sigma(z)=(r+s-z_1,z_2,\ldots,z_N)\). Then,
    \begin{equation}\label{eq:decreasing}
 \begin{aligned}
 (I_\alpha*Q^p)(re_1)&-(I_\alpha*Q^p)(se_1)\\
  &=\int_H [I_\alpha(re_1-z)-I_\alpha(se_1-z)] [Q(z)^p-Q(\sigma(z))^p]\,dz>0.
 \end{aligned}
   \end{equation}
  Indeed, on \(H\),
  \[
\begin{aligned}
 |se_1-z|^2-|re_1-z|^2&=(s-r)(r+s-2z_1)>0,\\
  |\sigma(z)|^2-|z|^2&=(r+s)(r+s-2z_1)>0.
\end{aligned}
 \]
 Since \(p\ge2\), the strict decrease of \(Q\) and \eqref{eq:decreasing}
 imply
  \begin{equation}\label{eq:decreasing-local-potential}
  r\longmapsto (I_\alpha*Q^p)(r)Q(r)^{p-2} \quad\text{is strictly decreasing on }(0,\infty).
    \end{equation}

 We first prove \eqref{eq:positive-unweighted-mean}. Since \(\eta\) is nondecreasing,
 nonconstant, and has zero \(\nu\)-mean, there exists \(r_0>0\) such that \(\eta(r)\le0\)
 for \(r\le r_0\) and \(\eta(r)\ge0\) for \(r\ge r_0\). Together with
 \eqref{eq:decreasing-local-potential}, this gives
\[
\begin{aligned}
   0=\int_{\R^N}\eta\,d\nu(x) ={}&\int_{\R^N}\eta Q^2 \bigl[(I_\alpha*Q^p)Q^{p-2} -(I_\alpha*Q^p)(r_0)Q(r_0)^{p-2}\bigr]\,dx\\
   &+(I_\alpha*Q^p)(r_0)Q(r_0)^{p-2} \int_{\R^N}\eta Q^2\,dx\\
 <{}&(I_\alpha*Q^p)(r_0)Q(r_0)^{p-2} \int_{\R^N}\eta Q^2\,dx.
    \end{aligned}
   \]
 This proves \eqref{eq:positive-unweighted-mean}.

 We next prove \eqref{eq:negative-dilation-weight}. We regularize the kernel before differentiating
  and symmetrizing. For \(\varepsilon>0\), replace \(I_\alpha(z)\) by
  \(I_{\alpha,\varepsilon}(z)=c_{N,\alpha}(|z|^2+\varepsilon^2)^{-(N-\alpha)/2}\). Since \(\nabla Q^p\in
  L^1(\R^N)\cap L^\infty(\R^N)\), for every \(\delta>0\),
\begin{equation}
  \label{eq:unif-conv-Riesz}
   \begin{aligned}
 &\|\nabla(I_{\alpha,\varepsilon}*Q^p) -\nabla(I_\alpha*Q^p)\|_{L^\infty(\R^N)}\\
 &\quad\le \|\nabla Q^p\|_{L^\infty(\R^N)} \int_{|z|<\delta}I_\alpha(z)\,dz + \|\nabla Q^p\|_{L^1(\R^N)} \sup_{|z|\ge\delta}
   |I_{\alpha,\varepsilon}(z)-I_\alpha(z)|\\
&\quad\le C\delta^\alpha+ C\varepsilon^2\delta^{\alpha-N-2}\to 0, \text{ first letting }
\varepsilon\to0 \text{ and then } \delta\to0.
    \end{aligned}
  \end{equation}
  By Lemma~\ref{lem:decay} and \eqref{eq:unif-conv-Riesz},
  \begin{equation}\label{eq:lim-varepsilon}
  \begin{aligned}
  &\left| \int_{\R^N}\eta Q^p x\cdot \left[ \nabla(I_{\alpha,\varepsilon}*Q^p) -\nabla(I_\alpha*Q^p) \right]\,dx \right|\\
  &\quad\le \|Q\eta\|_{L^2(\R^N)} \bigl\||x|Q^{p-1}\bigr\|_{L^2(\R^N)} \bigl\| \nabla(I_{\alpha,\varepsilon}*Q^p)
  -\nabla(I_\alpha*Q^p) \bigr\|_{L^\infty(\R^N)} \to0.
  \end{aligned}
\end{equation}
 Using
   \[
   \nabla I_{\alpha,\varepsilon}(x-y) = -(N-\alpha)c_{N,\alpha} \frac{x-y} {(|x-y|^2+\varepsilon^2)^{(N-\alpha+2)/2}}
    \]
  and symmetry in \(x\) and \(y\), we obtain
  \begin{equation}\label{eq:regularity-identity}
   \begin{aligned}
 &\int_{\R^N}\eta Q^p x\cdot\nabla(I_{\alpha,\varepsilon}*Q^p)\,dx\\
  &\quad= -\frac{N-\alpha}{4} \iint_{\R^N\times\R^N} c_{N,\alpha}Q(x)^pQ(y)^p
  \frac{[\eta(|x|)+\eta(|y|)]|x-y|^2} {(|x-y|^2+\varepsilon^2)^{(N-\alpha+2)/2}} \,dx\,dy\\
&\qquad -\frac{N-\alpha}{4} \iint_{\R^N\times\R^N} c_{N,\alpha}Q(x)^pQ(y)^p \frac{ [\eta(|x|)-\eta(|y|)]
(|x|^2-|y|^2) } {(|x-y|^2+\varepsilon^2)^{(N-\alpha+2)/2}} \,dx\,dy\\
  &\quad=:I+II.
\end{aligned}
 \end{equation}
    Using symmetry in \(x\) and \(y\) together with
  \[
  \frac{|x-y|^2} {(|x-y|^2+\varepsilon^2)^{(N-\alpha+2)/2}} = \frac{1} {(|x-y|^2+\varepsilon^2)^{(N-\alpha)/2}} -
    \frac{\varepsilon^2} {(|x-y|^2+\varepsilon^2)^{(N-\alpha+2)/2}},
   \]
  we have
 \[
\begin{aligned}
  I&= -\frac{N-\alpha}{2} \int_{\R^N} \eta Q^p(I_{\alpha,\varepsilon}*Q^p)\,dx\\
   &\quad+ \frac{(N-\alpha)c_{N,\alpha}}{2} \int_{\R^N}\eta(|x|)Q(x)^p \int_{\R^N}
   \frac{\varepsilon^2Q(y)^p} {(|x-y|^2+\varepsilon^2)^{(N-\alpha+2)/2}} \,dy\,dx\\
 &\to -\frac{N-\alpha}{2} \int_{\R^N}\eta Q^p(I_\alpha*Q^p)\,dx = -\frac{N-\alpha}{2}\int_{\R^N}\eta\,d\nu =0.
 \end{aligned}
    \]
    Indeed, \(I_{\alpha,\varepsilon}\le I_\alpha\) and \(I_{\alpha,\varepsilon}\to
    I_\alpha\) pointwise give, by dominated convergence,
 \[
    \int_{\R^N}\eta Q^p(I_{\alpha,\varepsilon}*Q^p)\,dx \to \int_{\R^N}\eta Q^p(I_\alpha*Q^p)\,dx.
  \]
   And,
\[
 \begin{aligned}
\sup_{x\in\R^N}& \int_{\R^N} \frac{\varepsilon^2Q(y)^p} {(|x-y|^2+\varepsilon^2)^{(N-\alpha+2)/2}} \,dy\\
  &\quad\le \|Q\|_{L^\infty(\R^N)}^p\varepsilon^2 \int_{|z|<1} \frac{dz} {(|z|^2+\varepsilon^2)^{(N-\alpha+2)/2}}
    + \varepsilon^2\|Q^p\|_{L^1(\R^N)}\\
   &\quad\le C\varepsilon^\alpha \int_0^{1/\varepsilon} \frac{s^{N-1}} {(1+s^2)^{(N-\alpha+2)/2}} \,ds + C\varepsilon^2 =o(1).
   \end{aligned}
 \]

  On the other hand, since \(\eta\) is increasing,
  \([\eta(|x|)-\eta(|y|)](|x|^2-|y|^2)\ge0\), and therefore monotone convergence gives
\[
 II\downarrow -\frac{N-\alpha}{4} \iint_{\R^N\times\R^N} I_\alpha(x-y)Q(x)^pQ(y)^p \frac{ [\eta(|x|)-\eta(|y|)] (|x|^2-|y|^2) }
 {|x-y|^2} \,dx\,dy.
\]
   Letting \(\varepsilon\to0\) in \eqref{eq:regularity-identity} together with
    \eqref{eq:lim-varepsilon} gives
  \begin{equation}\label{eq:eta-Qp<0}
   \begin{aligned}
  &\int_{\R^N}\eta Q^p x\cdot\nabla(I_\alpha*Q^p)\,dx\\
  &=-\frac{N-\alpha}{4}\iint_{\R^N\times\R^N}
  I_\alpha(x-y)Q(x)^pQ(y)^p\frac{[\eta(|x|)-\eta(|y|)](|x|^2-|y|^2)}{|x-y|^2}dxdy<0.
\end{aligned}
   \end{equation}
   By~\eqref{eq:A0-ground-form}, Lemma~\ref{lem:decay}, and the Cauchy--Schwarz inequality,
  \[
   \int_{\R^N}|x||\eta'(|x|)|\,d\nu \le \left(\int_{\R^N}Q^2|\nabla\eta|^2\,dx\right)^{1/2}
   \bigl\||x|(I_\alpha*Q^p)Q^{p-1}\bigr\|_{L^2(\R^N)} <\infty.
    \]
 Moreover, the decay of \(Q\) and \(\nabla Q\), together with \(\nabla(I_\alpha*Q^p)=I_\alpha*\nabla
 Q^p\in L^\infty(\R^N)\), gives \(|x|Q\nabla\bigl((I_\alpha*Q^p)Q^{p-2}\bigr)\in L^2(\R^N)\). Thus all terms below
 are integrable. Integration by parts, \(\int_{\R^N}\eta\,d\nu=0\), and \eqref{eq:eta-Qp<0} therefore give
 \[
  \begin{aligned}
 &\int_{\R^N} \eta Q^2x\cdot\nabla\bigl((I_\alpha*Q^p)Q^{p-2}\bigr)\,dx\\
&= \int_{\R^N}\eta Q^p x\cdot\nabla(I_\alpha*Q^p)\,dx + \frac{p-2}{p}
 \int_{\R^N}\eta(I_\alpha*Q^p)x\cdot\nabla(Q^p)\,dx\\
   &= \int_{\R^N}\eta Q^p x\cdot\nabla(I_\alpha*Q^p)\,dx -\frac{p-2}{p} \int_{\R^N}\eta Q^p x\cdot\nabla(I_\alpha*Q^p)\,dx -\frac{p-2}{p}
   \int_{\R^N}|x|\eta'(|x|)\,d\nu\\
    &= \frac{2}{p} \int_{\R^N}\eta Q^p x\cdot\nabla(I_\alpha*Q^p)\,dx -\frac{p-2}{p} \int_{\R^N}|x|\eta'(|x|)\,d\nu <0.
 \end{aligned}
 \]
This proves~\eqref{eq:negative-dilation-weight}.
\end{proof}

\subsection{Radial coercivity and nondegeneracy}

 \begin{proposition}
   \label{prop:radial-coercivity}
 There exists \(c_Q>0\) such that
  \begin{equation}
   \label{eq:coercive-LQ}
   \langle h,L_Qh\rangle\ge c_Q\|h\|_{H^1(\R^N)}^2 \qquad(h\in T).
 \end{equation}
   Moreover,
   \[
    \operatorname{ind}_{\mathrm{rad}}(L_Q)=1,\qquad \ker L_Q\cap H^1_{\mathrm{rad}}(\R^N)=\{0\}.
\]
\end{proposition}

\begin{proof}
    By Lemma~\ref{lem:ground-state-coercivity},
\begin{equation}\label{eq:constrained-maximum}
0<\tau=\sup_{h\in T\setminus\{0\}} \frac{\langle h,(A_0-L_Q)h\rangle}{\langle h,A_0h\rangle}
  \end{equation}
  is attained. Positivity follows by choosing a nonconstant radial \(\eta\), with \(\eta'\ge0\) compactly
    supported in \((0,\infty)\) and \(\int_{\R^N}\eta\,d\nu(x)=0\), and applying \eqref{eq:covariance-form}
    and \eqref{eq:positive-covariance}. Let \(h=Q\eta\) attain \eqref{eq:constrained-maximum}. The function
    \(\widetilde\eta\) in \eqref{eq:def-title-eta} satisfies \(\langle
    Q\widetilde\eta,A_0(Q\widetilde\eta)\rangle=\langle Q\eta,A_0(Q\eta)\rangle\), whereas
    \eqref{eq:sign-increases} implies \(Q\widetilde\eta\) increases \eqref{eq:constrained-maximum} if
    \(\eta'\) changes sign. We may therefore take
\[
  \eta'\ge0,\qquad \eta\text{ nonconstant},\qquad \int_{\R^N}\eta\,d\nu(x)=0.
  \]
  The Euler--Lagrange equation associated with \eqref{eq:constrained-maximum} is
   \begin{equation}\label{eq:identity-var}
  \langle\phi,(A_0-L_Q)h\rangle =\tau\langle\phi,A_0h\rangle\qquad(\phi\in T).
 \end{equation}

  We next extend \eqref{eq:identity-var} to \(\phi\in H^1_{\mathrm{rad}}(\R^N)\).
  Firstly,
  \begin{equation}\label{eq:decomp-phi}
   \phi=\left(\phi- \frac{\int_{\R^N}(\phi/Q)\,d\nu(x)}{\nu(\R^N)}Q\right)
   +\frac{\int_{\R^N}(\phi/Q)\,d\nu(x)}{\nu(\R^N)}Q \in T\oplus\operatorname{span}\{Q\}.
\end{equation}
  Since \(A_0Q=0\) and \(\langle Q,(A_0-L_Q)h\rangle
  =2(p-1)\int_{\R^N}\eta\,d\nu(x)=0\), so
 \begin{equation}\label{eq:generalized-eigenfunction}
    \tau A_0h=(A_0-L_Q)h \qquad\text{in }H^{-1}_{\mathrm{rad}}(\R^N).
  \end{equation}

  We claim that \(\tau<1\). For \(\omega>0\), set
\[
  Q_\omega(x) = \omega^{(\alpha+2)/(4(p-1))} Q(\sqrt{\omega}\,x).
   \]
Then
  \[
  -\Delta Q_\omega+\omega Q_\omega = (I_\alpha*Q_\omega^p)Q_\omega^{p-1}.
   \]
  Differentiating at \(\omega=1\), we obtain in \(H^{-1}(\R^N)\)
   \begin{equation}\label{eq:frequency}
  L_Q\left( \frac{\alpha+2}{4(p-1)}Q +\frac12x\cdot\nabla Q \right) = -Q.
   \end{equation}

 Testing~\eqref{eq:generalized-eigenfunction} against
  \(\frac{\alpha+2}{4(p-1)}Q+\frac12x\cdot\nabla Q\), and
  using~\eqref{eq:frequency} and self-adjointness, we obtain
\[
  (\tau-1) \left\langle h, A_0\left( \frac{\alpha+2}{4(p-1)}Q+\frac12x\cdot\nabla Q \right) \right\rangle = \int_{\R^N}hQ\,dx.
\]
 Since \(A_0Q=0\) and
\[
 A_0(x\cdot\nabla Q) = Q\,x\cdot\nabla\bigl((I_\alpha*Q^p)Q^{p-2}\bigr) +
 2\bigl((I_\alpha*Q^p)Q^{p-2}-1\bigr)Q,
  \]
together with \(h=Q\eta\) and \(\int_{\R^N}\eta\,d\nu=0\) give
  \[
    \left\langle h, A_0\left( \frac{\alpha+2}{4(p-1)}Q+\frac12x\cdot\nabla Q \right) \right\rangle = \frac12 \int_{\R^N}\eta
    Q^2 x\cdot\nabla\bigl((I_\alpha*Q^p)Q^{p-2}\bigr)\,dx - \int_{\R^N}Q^2\eta\,dx.
\]
    Hence
  \begin{equation}\label{eq:dilation-eigenvalue}
 \int_{\R^N}Q^2\eta\,dx = \frac12\left(1-\frac1\tau\right) \int_{\R^N}\eta Q^2 x\cdot\nabla\bigl((I_\alpha*Q^p)Q^{p-2}\bigr)\,dx.
    \end{equation}
   Together with Lemma~\ref{lem:radial-virial-signs}, \(0<\tau<1\).

  By \eqref{eq:constrained-maximum},
\begin{equation}\label{eq:strict-radial-stability}
   \langle h,L_Qh\rangle \ge(1-\tau)\langle h,A_0h\rangle>0 \qquad(h\in T\setminus\{0\}).
\end{equation}
Combining \eqref{eq:strict-radial-stability} with
Lemma~\ref{lem:ground-state-coercivity} gives
\eqref{eq:coercive-LQ}. By \eqref{eq:decomp-phi}, one has
\(H^1_{\mathrm{rad}}(\R^N)=T\oplus\operatorname{span}\{Q\}\).
Moreover,
  \[\langle Q,L_QQ\rangle
  =-2(p-1)\int_{\R^N}Q^p(I_\alpha*Q^p)\,dx<0.
  \]
    Therefore,
 \[
    \operatorname{ind}_{\mathrm{rad}}(L_Q)=1.
  \]
 If \(g\in\ker L_Q\cap H^1_{\mathrm{rad}}(\R^N)\), then
   \[
0 = \langle g,L_QQ\rangle = -2(p-1) \int_{\R^N}(I_\alpha*Q^p)Q^{p-1}g\,dx.
  \]
 Thus \(g\in T\), and~\eqref{eq:strict-radial-stability} implies
    \[
    \ker L_Q\cap H^1_{\mathrm{rad}}(\R^N)=\{0\}.
   \]
  This completes the proof.
  \end{proof}

\subsection{Nonradial sectors}

    Proposition~\ref{prop:radial-coercivity} treats the radial sector. Let
    \(N\ge2\), and let \(\{Y_{\ell,m}\}\) be a real orthonormal
    basis of spherical harmonics on \(\mathbb S^{N-1}\). For
    every \(h\in H^1(\mathbb R^N)\), polar coordinates show
    that \(h(r\,\cdot)\in L^2(\mathbb S^{N-1})\) for almost every \(r>0\).
    Consequently,
\[
 h(r\theta)=\sum_{\ell,m}f_{\ell,m}(r)Y_{\ell,m}(\theta), \qquad f_{\ell,m}(r)=\int_{\mathbb S^{N-1}}
 h(r\omega)Y_{\ell,m}(\omega)\,d\omega.
  \]
 Orthogonality of the spherical harmonics gives
 \begin{equation}
  \label{eq:orth-1}
   \begin{aligned}
    &\int_{\R^N}\left( |\nabla h|^2+h^2 -(p-1)(I_\alpha*Q^p)Q^{p-2}h^2 \right)\,dx\\
    &\quad= \sum_{\ell,m}\int_0^\infty \big[ |f_{\ell,m}'(r)|^2 +\frac{\ell(\ell+N-2)}{r^2}|f_{\ell,m}(r)|^2 \\
 &\qquad \qquad+\bigl(1-(p-1)(I_\alpha*Q^p)(r)Q(r)^{p-2}\bigr) |f_{\ell,m}(r)|^2 \big]r^{N-1}\,dr.
\end{aligned}
  \end{equation}

    Moreover, for fixed \(r,s>0\) with \(r\ne s\),
 \[
   I_\alpha(r\theta-s\omega) = c_{N,\alpha} (r^2+s^2-2rs\,\theta\cdot\omega)^{-(N-\alpha)/2}
 \]
    is a bounded function of \(\theta\cdot\omega\). The Funk--Hecke
    formula \cite[Theorem~1.2.9]{DX13} therefore implies 
 \[
\theta\longmapsto \int_{\mathbb S^{N-1}} I_\alpha(r\theta-s\omega)Y_{\ell,m}(\omega)\,d\omega
\]
  is a scalar multiple of \(Y_{\ell,m}(\theta)\). By Fubini's theorem and orthogonality of
  the spherical harmonics, for \((\ell,m)\ne(k,j)\), we obtain
 \begin{equation}
  \label{eq:orth-2}
   \begin{aligned}
   \int_{\R^N} & \bigl(I_\alpha*(Q^{p-1}f_{\ell,m}Y_{\ell,m})\bigr) Q^{p-1}f_{k,j}Y_{k,j}\,dx\\
&= \int_0^\infty\int_0^\infty Q(s)^{p-1}f_{\ell,m}(s) Q(r)^{p-1}f_{k,j}(r)\,s^{N-1}r^{N-1}\\
 &\qquad{}\cdot \int_{\mathbb S^{N-1}} \left[ \int_{\mathbb S^{N-1}} I_\alpha(r\theta-s\omega) Y_{\ell,m}(\omega)\,d\omega \right]
 Y_{k,j}(\theta)\,d\theta\,ds\,dr\\
  &=0.
  \end{aligned}
   \end{equation}
Here the angular cancellation holds for \(r\ne s\); the diagonal \(r=s\) has
  measure zero in the radial integration.

 Therefore, by \eqref{eq:orth-1} and \eqref{eq:orth-2},
 \begin{equation}
\label{eq:decomp-LQ-harmonic}
 \langle L_Qh,h\rangle =\sum_{\ell,m} \langle L_Q(f_{\ell,m}Y_{\ell,m}), f_{\ell,m}Y_{\ell,m}\rangle.
  \end{equation}
Differentiating~\eqref{eq:equation} in the distributional sense gives
   \[
 L_Q(\partial_jQ)=0, \qquad \partial_jQ(x)=Q'(|x|)\frac{x_j}{|x|}.
\]
Hence the translation modes lie in degree one. For \(x=r\theta\),
   \begin{equation}\label{eq:lap-sep}
   -\Delta(f(r)Y_{\ell,m}(\theta)) =\left(-f''-\frac{N-1}{r}f'
   +\frac{\ell(\ell+N-2)}{r^2}f\right)Y_{\ell,m}(\theta).
\end{equation}
  For \(f\in C_c^\infty((0,\infty))\), \cite[Definition~2.1 and Proposition~4.3]{BM24}, applied to
  the radial profile \(u(r)=r^{-\ell}f(r)\), gives
    \[
 e^{t\Delta}(fY_{\ell,m})(r\theta) = Y_{\ell,m}(\theta) \int_0^\infty p_\ell(t;r,s)f(s)s^{N-1}\,ds,
 \]
   where
  \[
p_\ell(t;r,s) = \frac{(rs)^{-(N-2)/2}}{2t} e^{-(r^2+s^2)/(4t)} I_{\ell+(N-2)/2}\left(\frac{rs}{2t}\right),
 \]
  and \(I_\nu\) denotes the modified Bessel function of the first kind.
  Moreover, the heat-semigroup representation of the standard Riesz potential
  \cite[Theorem~5]{Stinga19} gives
    \begin{equation}
  \label{eq:id-1}
  \begin{aligned}
I_\alpha*(fY_{\ell,m})(r\theta) &= \frac{1}{\Gamma(\alpha/2)} \int_0^\infty t^{\alpha/2-1}e^{t\Delta}(fY_{\ell,m})(r\theta)\,dt
\\
  &= Y_{\ell,m}(\theta) \int_0^\infty \kappa_\ell(r,s)f(s)s^{N-1}\,ds,
   \end{aligned}
 \end{equation}
 where
\[
   \kappa_\ell(r,s) = \frac{1}{\Gamma(\alpha/2)} \int_0^\infty t^{\alpha/2-1}p_\ell(t;r,s)\,dt.
  \]
    For fixed \(z>0\), \(I_\nu(z)\) is positive and strictly decreasing with
      respect to \(\nu>0\); see~\cite[Section~10.37]{OLBC10}.
      Hence
   \begin{equation}\label{eq:angular-kernels}
0<\kappa_{\ell+1}(r,s)<\kappa_\ell(r,s)
 \qquad(r,s>0,\ r\ne s).
   \end{equation}

 By Lemma~\ref{lem:decay}, the operator
 \(Q^{p-1}\bigl(I_\alpha*(Q^{p-1}\,\cdot)\bigr)\) is bounded on \(L^2(\R^N)\).
 Applying \eqref{eq:id-1} to \(Q^{p-1}f\), first for \(f\in
 C_c^\infty((0,\infty))\) and then by density, gives \eqref{eq:id-1} for
 every \(f\in L^2((0,\infty),r^{N-1}\,dr)\). Consequently, by
 \eqref{eq:linearization} and \eqref{eq:lap-sep}, \(L_Q(fY_{\ell,m})=(L_\ell
 f)Y_{\ell,m}\), where
\begin{equation}\label{eq:def-L}
 \begin{aligned}
 L_\ell f(r)={}&-f''(r)-\frac{N-1}{r}f'(r) +\frac{\ell(\ell+N-2)}{r^2}f(r)+f(r)\\
&-(p-1)(I_\alpha*Q^p)(r)Q(r)^{p-2}f(r)\\
&-pQ(r)^{p-1}\int_0^\infty \kappa_\ell(r,s)Q(s)^{p-1}f(s)s^{N-1}\,ds.
   \end{aligned}
\end{equation}

\begin{proposition}
 \label{prop:nonradial-sectors}
 Assume \(N\ge2\), and write
   \[
  q_\ell(f)=\langle f,L_\ell f\rangle_{L^2((0,\infty),r^{N-1}\,dr)}.
  \]
Then
 \[
   q_1\ge0,\qquad \ker L_1=\operatorname{span}\{Q'\},
 \]
and \(q_\ell(f)>0\) for every \(\ell\ge2\) and every nonzero \(f\). For \(N=1\),
   \[
 \langle h,L_Qh\rangle\ge0 \quad(h\in H^1_{\mathrm{odd}}(\R)),\qquad \ker L_Q\cap
   H^1_{\mathrm{odd}}(\R)=\operatorname{span}\{Q'\}.
    \]
\end{proposition}

 \begin{proof}
 \emph{Dimensions \(N\ge2\).} By Lemma~\ref{lem:decay}, \(\sigma_{\mathrm{ess}}(L_\ell)=[1,\infty)\). A direct
   calculation gives
  \begin{equation}\label{eq:angular-modulus}
    \begin{aligned}
  &q_\ell(|f|)-q_\ell(f)\\
 &= -p\iint_{(0,\infty)^2} Q(r)^{p-1}Q(s)^{p-1}\kappa_\ell(r,s) \Bigl[ |f(r)||f(s)|-f(r)f(s) \Bigr]
  (rs)^{N-1}\,dr\,ds\\
  &= -4p\iint_{(0,\infty)^2} Q(r)^{p-1}f_+(r)\kappa_\ell(r,s) Q(s)^{p-1}f_-(s) r^{N-1}s^{N-1}\,dr\,ds\le 0,
   \end{aligned}
   \end{equation}
with strict inequality if \(f\) changes sign.

We claim that \(q_1\ge0\). Suppose, to the
 contrary, that \(q_1(f)<0\) for some \(f\ne0\). Then
\[
  \lambda_1:=\inf\sigma(L_1) \le \frac{q_1(f)} {\|f\|_{L^2((0,\infty),r^{N-1}\,dr)}^2} <0.
\]
    Since \(\inf\sigma_{\mathrm{ess}}(L_1)=1\), \(\lambda_1\) is a discrete eigenvalue
      of \(L_1\). Thus there exists a real \(f\ne0\) such that \(L_1f=\lambda_1f\). If \(f\)
    changes sign, then \eqref{eq:angular-modulus} would give
\[
 q_1(|f|)<q_1(f), \qquad \||f|\|_{L^2((0,\infty),r^{N-1}dr)} = \|f\|_{L^2((0,\infty),r^{N-1}dr)},
   \]
 contradicting the minimality of \(\lambda_1\). Hence \(f\) has
  one sign, and we may assume \(f\ge0\), \(f\not\equiv0\). Since \(Q'(r)<0\)
  for every \(r>0\), \(\int_0^\infty f(r)Q'(r)r^{N-1}\,dr<0\). On the other
  hand, \(L_1f=\lambda_1f\), \(L_1Q'=0\), and self-adjointness
  of \(L_1\) give \(\lambda_1 \int_0^\infty f(r)Q'(r)r^{N-1}\,dr=0\). Since
 \(\lambda_1<0\), this is a contradiction. Therefore \(q_1\ge0\).

    Let \(f\in\ker L_1\). By \(q_1\ge0\) and \eqref{eq:angular-modulus}, we have \(0\le
      q_1(|f|)\le q_1(f)=0\). The strict part of~\eqref{eq:angular-modulus} shows that every
      nonzero \(f\in\ker L_1\) has a fixed sign. Since \(L_1Q'=0\), \(-Q'>0\), the kernel is
      one-dimensional. Indeed, otherwise one could choose \(0\ne f\in\ker L_1\)
      such that \(\int_0^\infty fQ'r^{N-1}\,dr=0\), which is impossible for a
      nonzero one-sign \(f\). Hence
   \[
 \ker L_1=\operatorname{span}\{Q'\}.
    \]
  For \(\ell\ge2\), first apply \eqref{eq:angular-modulus} to replace \(f\) by
    \(|f|\). Then \eqref{eq:angular-kernels}, \eqref{eq:def-L}, and \(q_1\ge0\)
  give, for \(f\ne0\),
\[
  \begin{aligned}
    q_\ell(f) &\ge q_\ell(|f|)\\
   &=q_1(|f|) +\bigl[\ell(\ell+N-2)-(N-1)\bigr] \int_0^\infty |f(r)|^2r^{N-3}\,dr\\
    &\quad+p\iint_{(0,\infty)^2} Q(r)^{p-1}Q(s)^{p-1} \bigl[\kappa_1(r,s)-\kappa_\ell(r,s)\bigr]
    |f(r)||f(s)|r^{N-1}s^{N-1}\,dr\,ds\\
&>0.
 \end{aligned}
  \]

\emph{Dimension \(N=1\).} It remains to analyze the odd sector. For
  \(h(x)=\operatorname{sgn}(x)f(|x|)\), \(f\in H^1_0(0,\infty)\),
\[
 (I_\alpha*(Q^{p-1}h))(r) =c_{1,\alpha}\int_0^\infty
 [|r-s|^{\alpha-1}-(r+s)^{\alpha-1}]Q(s)^{p-1}f(s)\,ds.
 \]
  Hence, by \eqref{eq:linearization} and \(1>\alpha>0\),
 \begin{equation}\label{eq:1D-identify}
\begin{aligned}
  &\left\langle\operatorname{sgn}(x)|f(|x|)|, L_Q\bigl(\operatorname{sgn}(x)|f(|x|)|\bigr)\right\rangle -\langle h,L_Qh\rangle\\
  &=-8pc_{1,\alpha}\iint_{(0,\infty)^2} [|r-s|^{\alpha-1}-(r+s)^{\alpha-1}]
    Q(r)^{p-1}Q(s)^{p-1}f_+(r)f_-(s)\,dr\,ds\le0,
  \end{aligned}
  \end{equation}
strictly if \(f\) changes sign. Suppose that the odd sector had a negative
direction. The same restriction argument gives essential spectrum
 \([1,\infty)\) in this sector. Hence its lowest spectral value
 \(\lambda<0\) is an eigenvalue. By~\eqref{eq:1D-identify}, a corresponding
 eigenfunction \(h\) may be chosen to have one sign on \((0,\infty)\); after
changing sign, let \(h\ge0\) there. Since
\[
  L_QQ'=0, \qquad Q'(r)<0 \quad(r>0),
  \]
 self-adjointness gives \(0=\langle L_Qh,Q'\rangle=\lambda\langle h,Q'\rangle\),
  whereas \(\langle h,Q'\rangle=2\int_0^\infty h(r)Q'(r)\,dr<0\), a contradiction.
  Hence
  \begin{equation}\label{eq:odd-nonnegative}
 \langle h,L_Qh\rangle\ge0 \qquad \text{for every }h\in H^1_{\mathrm{odd}}(\R).
\end{equation}
By~\eqref{eq:1D-identify} and~\eqref{eq:odd-nonnegative}, every nonzero element of \(\ker
L_Q\cap H^1_{\mathrm{odd}}(\R)\) has one sign on \((0,\infty)\). Since \(Q'\in\ker
L_Q\cap H^1_{\mathrm{odd}}(\R)\), this kernel is one-dimensional: otherwise one
could choose \(0\ne h\) in the kernel with \(\langle h,Q'\rangle=0\), while, after
changing sign, \(\langle h,Q'\rangle=2\int_0^\infty h(r)Q'(r)\,dr<0\), a contradiction.
Therefore
\[
   \ker L_Q\cap H^1_{\mathrm{odd}}(\R) = \operatorname{span}\{Q'\}.
  \]
 This completes the proof.
\end{proof}

\subsection{Proof of Theorem~\ref{thm:index-one-nondegeneracy}}

\begin{proof}[Proof of Theorem~\ref{thm:index-one-nondegeneracy}]

   Suppose first that \(N\ge2\). By Propositions \ref{prop:radial-coercivity}
   and \ref{prop:nonradial-sectors},
 \[
\begin{array}{c|c|c}
\text{degree} & \text{negative directions} & \text{kernel}\\ \hline
 \ell=0 & \text{one} & \{0\}\\
 \ell=1 & \text{none} & \operatorname{span}\{Q'(r)Y_{1,m}(\theta)\}_m\\
 \ell\ge2 & \text{none} & \{0\}.
\end{array}
    \]
  Since the degree-one spherical harmonics span \(\{\theta_1,\ldots,\theta_N\}\) and
  \(Q'(|x|)\frac{x_j}{|x|}=\partial_jQ(x)\). Therefore, by \eqref{eq:decomp-LQ-harmonic}, we have
\[
 \operatorname{ind}(L_Q)=1,\qquad \ker L_Q =\operatorname{span}\{\partial_1Q,\ldots,\partial_NQ\}.
\]

 For \(N=1\), Proposition \ref{prop:radial-coercivity} gives one negative direction and no kernel in the even
 subspace. Proposition \ref{prop:nonradial-sectors} gives no negative direction and kernel
 \(\operatorname{span}\{Q'\}\) in the odd subspace. Thus
  \[
   \operatorname{ind}(L_Q)=1,\qquad \ker L_Q=\operatorname{span}\{Q'\}.
  \]
This completes the proof.
\end{proof}

\section{Uniqueness}\label{sec:uniqueness}

\subsection{Uniqueness near the endpoint}

We now pass from radial nondegeneracy to uniqueness. For fixed \(p>2\), we
first show that every positive solution approaches the same profile
as \(\alpha\uparrow N\). Local uniqueness near that profile provides a
starting point for continuation throughout
    \[
  \max\{0,(N-2)p-N\}<\alpha<N.
    \]
The case \(p=2\) requires a separate argument, in which we continue in \(p\) from each quadratic
profile.

  For the endpoint limit, we remove the Riesz normalization constant. With \(\gamma=N-\alpha\), the
  substitution \( U=c_{N,\alpha}^{1/(2p-2)}Q\) into \eqref{eq:equation} gives
    \begin{equation}\label{eq:unnormalized}
  -\Delta U+U=(|\cdot|^{-\gamma}*U^p)U^{p-1} \qquad\text{in }\R^N.
  \end{equation}
 For each such
  \(U\), write \(V_\gamma=|\cdot|^{-\gamma}*U^p\).

\begin{proposition}
 \label{lem:endpoint}
    Let \(N\ge1\) and \(2<p<2^*\), where \(2^*=2N/(N-2)\) for \(N\ge3\) and \(2^*=\infty\) for
    \(N=1,2\). There exists \(\gamma_*>0\) such that, for every \(0<\gamma<\gamma_*\),
    equation~\eqref{eq:unnormalized} has exactly one positive solution in \(H^1(\R^N)\), up to
    translations.
\end{proposition}

 To pass to the limit \(\gamma\downarrow0\), we need estimates uniform over all centered positive
 solutions. We first bound \(U(0)^{p-2}V_\gamma(0)\), which controls the coefficient in
 \eqref{eq:unnormalized} for \(U/U(0)\). This estimate will yield upper and lower bounds for \(U(0)\) and
 uniform exponential decay of \(U\) and \(\nabla U\). These bounds provide compactness in
 \(H^1(\R^N)\cap L^\infty(\R^N)\) and exclude a zero limit. Throughout this subsection,
 fix \(N\ge1\), \(2<p<2^*\), and
  \[
  0<\bar\gamma<\min\{N/2,\,2N-(N-2)p\}.
\]

\begin{lemma}
 \label{lem:endpoint-potential}
    There exist \(\gamma_0\in(0,\bar\gamma]\) and \(C>0\) such that every
    centered positive solution of~\eqref{eq:unnormalized} with
    \(0<\gamma\le\gamma_0\) satisfies
  \begin{equation}\label{eq:endpoint-potential-bound}
  U(0)^{p-2}V_\gamma(0)\le C.
\end{equation}
\end{lemma}

\begin{proof}
 Suppose the assertion is false. Then there exist \(\gamma_j\downarrow0\) and centered
  positive solutions \(U_j\) such that \(U_j(0)^{p-2}V_{\gamma_j}(0)\to\infty\). Set
\[
\varepsilon_j= \bigl(U_j(0)^{p-2}V_{\gamma_j}(0)\bigr)^{-1/2}, \qquad u_j(x)=\frac{U_j(\varepsilon_jx)}{U_j(0)}.
\]
   Then \(\varepsilon_j\to0\), \(0<u_j\le1\), \(u_j(0)=1\), and
    \begin{equation}
\label{eq:appro-uj}
    -\Delta u_j+\varepsilon_j^2u_j = \frac{\displaystyle \int_{\R^N}|x-y|^{-\gamma_j}u_j(y)^p\,dy} {\displaystyle
    \int_{\R^N}|y|^{-\gamma_j}u_j(y)^p\,dy} u_j^{p-1}.
  \end{equation}
  The right of \eqref{eq:appro-uj} at most one, by radial monotonicity of the potential. Hence \(\Delta
  u_j\) is uniformly bounded. The interior estimate~\cite[Theorem~9.11]{GT83} gives, for every \(R>0\) and
  \(q>N\),
   \[
  \|u_j\|_{W^{2,q}(B_R)} \le C_R\bigl(\|u_j\|_{L^q(B_{2R})} +\|\Delta u_j\|_{L^q(B_{2R})}\bigr)\le C_R.
  \]
   Morrey embedding and \(u_j(0)=1\) imply that \(u_j\ge1/2\) on a fixed ball \(B_{r_0}\). Therefore
\begin{equation}\label{eq:endpoint-denominator}
 \int_{\R^N}|y|^{-\gamma_j}u_j(y)^p\,dy \ge2^{-p}\int_{B_{r_0}}|y|^{-\gamma_j}\,dy\ge c_0>0.
\end{equation}

  Fix \(K>0\) and \(R>2K+1\). On \(B_R\), use \(u_j\le1\). On its complement, use
  \[
\left| \left(\frac{|x-y|}{|y|}\right)^{-\gamma_j}-1 \right| \le C\frac{\gamma_jK}{R}
 \qquad(|x|\le K,\ |y|>R).
\]
   Together with~\eqref{eq:endpoint-denominator}, these bounds give
  \[
 \begin{aligned}
    &\sup_{|x|\le K} \left| \frac{\displaystyle\int_{\R^N} |x-y|^{-\gamma_j}u_j(y)^p\,dy} {\displaystyle\int_{\R^N}
    |y|^{-\gamma_j}u_j(y)^p\,dy}-1 \right| \le \frac{2}{c_0}\int_{B_{R+K}} \bigl||z|^{-\gamma_j}-1\bigr|\,dz +C\frac{\gamma_jK}{R}\to0.
\end{aligned}
   \]
  Consequently, after extraction, \(u_j\to u\) in
   \(C^1_{\mathrm{loc}}(\R^N)\), where
\[
-\Delta u=u^{p-1},\qquad u>0,\qquad u(0)=1.
\]
   For \(N\ge3\), this contradicts the subcritical Liouville theorem~\cite[Theorem~1.1]{GS81};
   see also \cite[Theorem~2]{CL91}. For \(N=1,2\), radiality yields
   \[
 -r^{N-1}u'(r) =\int_0^r t^{N-1}u(t)^{p-1}\,dt\ge c>0 \qquad(r\ge1).
  \]
  Integration of \(-u'(r)\ge cr^{1-N}\) contradicts positivity. This proves~\eqref{eq:endpoint-potential-bound}.
\end{proof}

\begin{lemma}
   \label{lem:endpoint-uniform-bounds}
  After reducing \(\gamma_0\) if necessary, there exist \(c,C,\eta_0>0\) such that every centered
  positive solution with \(0<\gamma\le\gamma_0\) satisfies
  \begin{equation}\label{eq:endpoint-bounds}
  \begin{aligned}
 c\le U(0)\le C,\quad |U(x)|+|\nabla U(x)|\le Ce^{-\eta_0|x|}, \quad \|U\|_{W^{2,r}(\R^N)}\le C_r
   \quad(1<r<\infty).
  \end{aligned}
   \end{equation}
   The constants are independent of \(U\) and \(\gamma\).
\end{lemma}

\begin{proof}
Set \(w=U/U(0)\). By Lemma~\ref{lem:endpoint-potential} and radial monotonicity, \(\Delta w\) is uniformly
bounded. The interior estimate~\cite[Theorem~9.11]{GT83} gives a fixed \(r_0>0\) such that \(U\ge U(0)/2\)
on \(B_{r_0}\). It follows that
\[
\begin{aligned}
 V_\gamma(0) &=\int_{\R^N}|y|^{-\gamma}U(y)^p\,dy \ge2^{-p}U(0)^p\int_{B_{r_0}}|y|^{-\gamma}\,dy \ge
 cU(0)^p.
    \end{aligned}
\]
 Hence
\begin{equation}\label{eq:endpoint-upper-center}
 U(0)^{2p-2} \le C U(0)^{p-2}V_\gamma(0)\le C.
  \end{equation}
 By \eqref{eq:unnormalized} and \(V_\gamma'\le0\),
   \begin{equation}\label{eq:radial-energy}
  \frac{d}{dr}\left[ \frac12U'(r)^2-\frac12U(r)^2 +\frac1pV_\gamma(r)U(r)^p\right] =-\frac{N-1}{r}U'(r)^2
    +\frac1pV_\gamma'(r)U(r)^p\le0.
\end{equation}
 By Lemma~\ref{lem:decay}, \( \frac12U'(r)^2-\frac12U(r)^2
 +\frac1pV_\gamma(r)U(r)^p\to0\). It is therefore nonnegative, and
\[
 V_\gamma(0)U(0)^{p-2}\ge\frac p2.
  \]
Choose
\[
  1<\kappa<\frac p2,\qquad R_*>\sqrt{\frac{\lambda_1(B_1)}{\kappa-1}},\qquad \eta_0=\sqrt{1-\frac{2\kappa}{p}},
\]
    where \(\lambda_1(B_1)\) is the first Dirichlet eigenvalue of \(-\Delta\) on
     \(B_1\). If \(V_\gamma(R_*)U(R_*)^{p-2}\ge\kappa\), testing \eqref{eq:unnormalized} with
     \(\phi^2/U\), for \(\phi\in C_c^\infty(B_{R_*})\), gives
    \[
  \begin{aligned}
 \int_{B_{R_*}}|\nabla\phi|^2\,dx &=\int_{B_{R_*}}U^2 \left|\nabla\left(\frac\phi U\right)\right|^2\,dx +\int_{B_{R_*}}(V_\gamma
  U^{p-2}-1)\phi^2\,dx \ge(\kappa-1)\int_{B_{R_*}}\phi^2\,dx.
  \end{aligned}
    \]
 Taking the infimum of the Rayleigh quotient contradicts \(\lambda_1(B_{R_*})=R_*^{-2}\lambda_1(B_1)<\kappa-1\).
   Therefore \(V_\gamma(r)U(r)^{p-2}<\kappa\) for \(r\ge R_*\). By \eqref{eq:radial-energy} and \(U'<0\) imply
 \[
 -\frac{U'(r)}{U(r)} \ge\sqrt{1-\frac2pV_\gamma(r)U(r)^{p-2}} \ge\eta_0 \qquad(r\ge R_*).
 \]
   Integrating, and using \(U(R_*)\le U(0)\), yields
\begin{equation}\label{eq:endpoint-envelope}
 U(r)\le U(0)\min\{1,e^{-\eta_0(r-R_*)}\}.
  \end{equation}
    Consequently,
\[
\begin{aligned}
V_\gamma(0) &\le U(0)^p\int_{\R^N}|y|^{-\gamma} \min\{1,e^{-p\eta_0(|y|-R_*)}\}\,dy \le C U(0)^p,
\end{aligned}
\]
 where \(C\) is uniform for \(0<\gamma\le\gamma_0\). Thus
  \[
    \frac p2\le V_\gamma(0)U(0)^{p-2} \le C U(0)^{2p-2},
   \]
  which proves \(U(0)\ge c>0\).

Lemma~\ref{lem:endpoint-potential} also gives \(|\Delta U|\le CU\). By
  \eqref{eq:endpoint-upper-center} and \eqref{eq:endpoint-envelope}, give uniform
  exponential bounds on \(U\) and \(\Delta U\). Interior gradient estimates on balls
  of fixed radius give the same bound on \(|\nabla U|\). Finally, the global elliptic
  estimate for \(-\Delta+1\) yields
   \[
  \|U\|_{W^{2,r}(\R^N)} \le C_r\|(-\Delta+1)U\|_{L^r(\R^N)}\le C_r, \qquad 1<r<\infty.
  \]
This proves~\eqref{eq:endpoint-bounds}.
 \end{proof}

 \begin{proof}[Proof of Proposition~\ref{lem:endpoint}]
 Let \(\gamma_j\downarrow0\), with corresponding centered positive
 solutions \(U_j\). By Lemma~\ref{lem:endpoint-uniform-bounds} and local
  compactness, up to a subsequence,
   \begin{equation}
  \label{eq:uj-to-U}
    U_j\to U_* \quad\text{strongly in }H^1(\R^N)\cap L^\infty(\R^N), \qquad U_*\ne0.
  \end{equation}
By \eqref{eq:endpoint-bounds} and dominated convergence, gives
\(V_{\gamma_j}\to\int_{\R^N}U_*^p\,dx\) locally uniformly; this also can
see \cite[Proposition~2.8]{Seok19}. Since \(U_j\to U_*\) uniformly, we may pass
to the limit in \eqref{eq:unnormalized} and obtain
    \begin{equation}
\label{eq:qstar}
  -\Delta U_*+U_* =\left(\int_{\R^N}U_*^p\,dx\right)U_*^{p-1}, \qquad U_*>0.
 \end{equation}
 By uniqueness theorem~\cite{MS87,Kwong89,McLeod93}
    \begin{equation}\label{eq:endpoint-profile}
U_*= \left(\int_{\R^N}w^p\,dx\right)^{-1/(2p-2)}w,
   \end{equation}
 where \(w\) satisfies
    \begin{equation}
  \label{eq:limit-equation}
 -\Delta w+w=w^{p-1}, \qquad w>0, \qquad w\in H^1_{\mathrm{rad}}(\R^N).
  \end{equation}
  Thus every centered positive sequence has the same limit \(U_*\).

 We finish by proving uniqueness near \(U_*\). The argument is
 motivated by the perturbative uniqueness scheme in
 \cite[Lemma~5.3]{Seok19}. The linearized operator of \eqref{eq:limit-equation}
 is
    \[A=-\Delta+1-(p-1)w^{p-2}.\]
 By the classical nondegeneracy of \(w\)~\cite{Kwong89,CGNT07},
   \(A:H^2_{\mathrm{rad}}(\R^N)\to L^2_{\mathrm{rad}}(\R^N)\) is invertible. The
   linearized operator of \eqref{eq:qstar} at \(U_*\), using
   \eqref{eq:endpoint-profile}, is
 \[
  L_*h =Ah- \frac{p w^{p-1}}{\displaystyle\int_{\R^N}w^p\,dx} \int_{\R^N}w^{p-1}h\,dx.
  \]
   Since \(Aw=-(p-2)w^{p-1}\), any radial \(h\) with \(L_*h=0\) satisfies
   \[
h= -\frac{p w}{(p-2)\displaystyle\int_{\R^N}w^p\,dx} \int_{\R^N}w^{p-1}h\,dx.
 \]
 Testing with \(w^{p-1}\) gives \(h=0\). As \(L_*\) is a rank-one perturbation of \(A\), the
  Fredholm alternative yields
  \begin{equation}\label{eq:endpoint-radial-invertibility}
   L_*:H^2_{\mathrm{rad}}(\R^N) \to L^2_{\mathrm{rad}}(\R^N) \quad\text{is invertible}.
\end{equation}

 We next transfer this invertibility to profiles approaching \(U_*\). Let
 \(L_{\gamma,u}\) be the linearization of~\eqref{eq:unnormalized} at
 a nonnegative radial function \(u\):
\[
 L_{\gamma,u}h = -\Delta h+h -(p-1)(|\cdot|^{-\gamma}*u^p)u^{p-2}h -pu^{p-1} \bigl(|\cdot|^{-\gamma}*(u^{p-1}h)\bigr).
\]
   By \eqref{eq:endpoint-profile}, \(L_*\) can
     equivalently be written as
\[
   \begin{aligned}
    L_*h ={}& -\Delta h+h -(p-1) \left(\int_{\R^N}U_*^p\,dx\right) U_*^{p-2}h -pU_*^{p-1} \int_{\R^N}U_*^{p-1}h\,dx.
 \end{aligned}
\]
   By the Cauchy--Schwarz inequality,
   \begin{equation}\label{eq:endpoint-operator}
 \begin{aligned}
 &\|L_{\gamma_j,u_j}-L_*\|_{L^2(\R^N)\to L^2(\R^N)} \\
  &\quad\le (p-1) \left\| (|\cdot|^{-\gamma_j}*u_j^p)u_j^{p-2}
    -\left(\int_{\R^N}U_*^p\,dx\right)U_*^{p-2} \right\|_{L^\infty(\R^N)} \\
 &\qquad +p\left( \iint_{\R^N\times\R^N} \left| u_j(x)^{p-1}|x-y|^{-\gamma_j}u_j(y)^{p-1}
   -U_*(x)^{p-1}U_*(y)^{p-1} \right|^2\,dx\,dy \right)^{1/2} \to0.
   \end{aligned}
    \end{equation}
    The first follows by \eqref{eq:uj-to-U} and \(0\le u_j(x)\le Ce^{-\eta_0|x|}\). The
      second follows by dominated convergence.

 Suppose that \(\gamma_j\downarrow0\) and that \(Q_j,R_j\) are distinct centered positive
 solutions. By \eqref{eq:uj-to-U}, \(Q_j,R_j\to U_*\) in \(H^1(\R^N)\cap
 L^\infty(\R^N)\). Set \(h_j=Q_j-R_j\) and \(u_{j,t}=R_j+th_j\), \(0\le t\le1\).
 Then, same as \eqref{eq:endpoint-operator},
  \begin{equation}\label{eq:unif-t}
  \sup_{0\le t\le1} \|L_{\gamma_j,u_{j,t}}-L_*\|_{L^2(\R^N)\to L^2(\R^N)}\to0.
  \end{equation}
 Define
\[
 F_\gamma(u) = -\Delta u+u-(|\cdot|^{-\gamma}*u^p)u^{p-1}.
\]
 Since
    \[
 F_{\gamma_j}(Q_j) = F_{\gamma_j}(R_j) =0, \qquad D_uF_{\gamma_j}(u_{j,t}) = L_{\gamma_j,u_{j,t}},
  \]
    the fundamental theorem of calculus gives
   \[
  \begin{aligned}
  0 &= F_{\gamma_j}(Q_j)-F_{\gamma_j}(R_j) = \int_0^1 L_{\gamma_j,u_{j,t}}h_j\,dt.
 \end{aligned}
  \]
 Therefore
\[
  L_*h_j = \int_0^1 (L_*-L_{\gamma_j,u_{j,t}})h_j\,dt.
 \]
  By~\eqref{eq:endpoint-radial-invertibility} and~\eqref{eq:unif-t},
    \[
  \|h_j\|_{L^2(\R^N)} \le C\|L_*h_j\|_{L^2(\R^N)} \le o(1)\|h_j\|_{L^2(\R^N)},
    \]
 a contradiction. Hence there is at most one centered positive solution for all
 sufficiently small \(\gamma>0\).

    Existence follows from \cite[Theorems~1 and 3]{MV13}.
     Hence, for all sufficiently small \(\gamma>0\),
     equation~\eqref{eq:unnormalized} has exactly one positive
     solution up to translations.
\end{proof}

\subsection{\texorpdfstring{Continuation in \(\alpha\)}{Continuation in a}}

  For fixed \(N\), the conditions in~\eqref{eq:range} with \(p>2\) are equivalent to
  \[
2<p<2^*, \qquad \max\{0,(N-2)p-N\}<\alpha<N.
    \]
  Fix an arbitrary \(p\in(2,2^*)\) and set
\[
J=J_p=\bigl(\max\{0,(N-2)p-N\},N\bigr).
 \]

   We prove uniqueness for every \(\alpha\in J\) by continuation in \(\alpha\). Proposition~\ref{lem:endpoint} gives
     uniqueness for \(\alpha\) sufficiently close to \(N\). At every positive solution, radial nondegeneracy
     gives a locally unique branch with respect to \(\alpha\). To extend uniqueness throughout \(J\), we also
     need compactness of the centered positive solutions when \(\alpha\) remains in a compact subinterval of \(J\).
     This follows from Li's a priori estimates \cite[Theorem~1.5]{Li22}. The next lemma combines these two properties
     to show that the number of centered positive solutions is constant on \(J\).

\begin{lemma}\label{lem:proper-continuation}
  Fix \(N\ge1\) and \(2<p<2^*\). The number of centered positive solutions \(Q\in H^1(\R^N)\)
  of~\eqref{eq:equation} is independent of \(\alpha\in J\).
\end{lemma}

\begin{proof}
  \emph{Step~1. Convergence of centered positive solutions as the
  parameter varies.} Fix \(\alpha_j\to\alpha_*\) with \(\alpha_*\in J\),
  and let \(Q_j\) be centered positive solutions corresponding to
  \(\alpha_j\). Since \(\alpha_*\) lies strictly inside the admissible
  range, the a priori estimates of~\cite[Theorem~1.5]{Li22} give
  \[
  \|Q_j\|_{W^{2,r}(\R^N)}\le C_r \qquad (1<r<\infty), \qquad |Q_j(x)|+|\nabla Q_j(x)| \le
   Ce^{-c|x|},
\]
   with constants independent of \(j\). Hence, after passing to a
    subsequence,
\[
  Q_j\to Q_* \qquad\text{strongly in } H^1(\R^N)\cap L^\infty(\R^N).
 \]
   Moreover,
\[
 \|Q_j\|_{H^1(\R^N)}^2 = \int_{\R^N}(I_{\alpha_j}*Q_j^p)Q_j^p\,dx \le C\|Q_j\|_{H^1(\R^N)}^{2p},
\]
 and therefore \(\|Q_j\|_{H^1(\R^N)}\ge c>0\). Passing to the limit in \eqref{eq:equation} gives
\[
  -\Delta Q_*+Q_* = (I_{\alpha_*}*Q_*^p)Q_*^{p-1}, \qquad Q_*>0.
   \]

\emph{Step~2. Local continuation of each centered positive solution.}
     Define
 \[
 F:J\times H^1_{\mathrm{rad}}(\R^N) \to H^{-1}_{\mathrm{rad}}(\R^N), \qquad F(\alpha,u) =-\Delta
 u+u-(I_\alpha*|u|^p)|u|^{p-2}u.
  \]
 By HLS and Sobolev inequalities, \(F\) is continuous in \((\alpha,u)\), continuously Fr\'echet
  differentiable with respect to \(u\), and that \(D_uF\) depends continuously on \((\alpha,u)\).

   Let \(Q\) be a centered positive solution at \(\alpha\in J\). The
   operator
  \[
    D_uF(\alpha,Q)=L_Q: H^1_{\mathrm{rad}}(\R^N) \to H^{-1}_{\mathrm{rad}}(\R^N)
\]
  is Fredholm of index zero. By Theorem~\ref{thm:index-one-nondegeneracy}, its radial kernel is
  trivial, so it is an isomorphism. The implicit function theorem therefore gives
  \(\varepsilon>0\), with \((\alpha-\varepsilon,\alpha+\varepsilon)\subset J\), and an open
  neighborhood \(\mathcal U\) of \(Q\) in \(H^1_{\mathrm{rad}}(\R^N)\) such that, for
  every \(\beta\in(\alpha-\varepsilon,\alpha+\varepsilon)\), the equation
  \(F(\beta,u)=0\) has exactly one solution \(u\in\mathcal U\). Denoting this solution by \(u_\beta\),
  the map \((\alpha-\varepsilon,\alpha+\varepsilon) \ni\beta\longmapsto
  u_\beta\in\mathcal U\) is continuous and satisfies \(u_\alpha=Q\).

    We claim that the local branch remains positive. Testing
\[
 -\Delta u_\beta+u_\beta=(I_\beta*|u_\beta|^p)|u_\beta|^{p-2}u_\beta
    \]
  with \(-u_\beta^-\) and applying the HLS and Sobolev inequalities gives
    \[
    \|u_\beta^-\|_{H^1(\R^N)}^2 =\int_{\R^N} (I_\beta*|u_\beta|^p)(u_\beta^-)^p\,dx \le C\|u_\beta\|_{H^1(\R^N)}^p
    \|u_\beta^-\|_{H^1(\R^N)}^p,
  \]
   uniformly for \(\beta\) near \(\alpha\). Since \(u_\beta\to Q>0\) in \(H^1(\R^N)\) and
   \(p>2\), this estimate forces \(u_\beta^-=0\) for \(\beta\) sufficiently close to
   \(\alpha\). The branch is nontrivial, so the strong maximum principle gives
   \(u_\beta>0\).

 \emph{Step~3. Constancy of the number of centered positive solutions.}
 There are only finitely many centered positive solutions at
 \(\alpha_0\). Indeed, otherwise we could choose a sequence of pairwise
 distinct such solutions \(Q_j\). Applying Step~1 with
 \(\alpha_j=\alpha_0\), we obtain, after passing to a subsequence, \(Q_j\to
 Q\quad\text{in }H^1(\R^N)\), where \(Q\) is a centered positive
 solution at \(\alpha_0\). By Step~2, \(Q\) is the only solution at
 \(\alpha_0\) in some neighborhood of \(Q\) in \(H^1_{\mathrm{rad}}(\R^N)\).
 Consequently, \(Q_j=Q\) for all sufficiently large \(j\), contradicting
 the pairwise distinctness of the \(Q_j\). Hence there are
 finitely many, denoted by \(Q_1,\ldots,Q_m\).

   By Step~2 and continuity of the local branches, we may choose pairwise disjoint
   open neighborhoods \(\mathcal U_i\) of \(Q_i\) in
   \(H^1_{\mathrm{rad}}(\R^N)\) and a common \(\varepsilon>0\), with
   \((\alpha_0-\varepsilon,\alpha_0+\varepsilon)\subset J\), such that, for every
   \(|\beta-\alpha_0|<\varepsilon\), each \(\mathcal U_i\) contains
   exactly one radial solution \(u_{i,\beta}\), which is positive. These
   solutions are distinct, so there are at least \(m\) centered positive
   solutions at each such \(\beta\).

  After reducing \(\varepsilon\), there are no further centered
   positive solutions. Otherwise, we could find
  \(\beta_j\to\alpha_0\) and centered positive solutions \(R_j\)
   at \(\beta_j\) satisfying \(R_j\notin\bigcup_{i=1}^m\mathcal U_i\).
  By Step~1, a subsequence converges in \(H^1(\R^N)\) to a centered positive
   solution at \(\alpha_0\), hence to some \(Q_i\). Since \(\mathcal U_i\) is
   open, \(R_j\in\mathcal U_i\) for all sufficiently large \(j\), a
   contradiction.

  Thus there are exactly \(m\) centered positive solutions for every
  \(\beta\) sufficiently close to \(\alpha_0\). Their number is
  therefore locally constant on \(J\), and hence constant because \(J\) is an
  interval.
\end{proof}

\begin{remark}
 The conclusion of Lemma~\ref{lem:proper-continuation} reflects the
   covering property of proper local homeomorphisms over an interval;
 see~\cite[Lemma~2]{Ho75}.
\end{remark}

By Lemma~\ref{lem:endpoint}, for \(\alpha\) sufficiently close to \(N\),
equation~\eqref{eq:equation} has exactly one centered positive solution. By
Lemma~\ref{lem:proper-continuation}, the number of centered positive solutions is
constant on \(J\). Hence this number is equal to one for every \(\alpha\in J\). Therefore the
centered positive solution is unique for every \(p>2\) in \eqref{eq:range}.

\subsection{The quadratic case}

 It remains to treat \(p=2\). We keep \(\alpha\) fixed and continue each
 centered positive solution to exponents \(p>2\).

\begin{lemma}\label{lem:quadratic-uniqueness}
 Let \(N\ge1\) and \(\max\{0,N-4\}<\alpha<N\). Equation~\eqref{eq:equation} with \(p=2\)
 has exactly one centered positive radial solution in \(H^1(\R^N)\).
\end{lemma}

\begin{proof}
    Existence follows from~\cite[Theorems~1 and~3]{MV13}. Fix a centered positive
      radial solution \(Q\) at \(p=2\). Choose \(\delta>0\) such that
      \((N-2)(2+\delta)<N+\alpha\). Set
 \[
X=H^1_{\mathrm{rad}}(\R^N)\cap L^\infty(\R^N), \qquad \|u\|_X=\|u\|_{H^1(\R^N)}+\|u\|_{L^\infty(\R^N)},
 \]
 and define
\[
F(p,u) = u-(-\Delta+1)^{-1} \left[(I_\alpha*|u|^p)|u|^{p-2}u\right], \qquad 2\le p\le2+\delta.
\]
    For fixed \(p\), the map \(F\) is continuously differentiable in \(u\). As
      \((p,u)\to(2,Q)\) in \([2,2+\delta]\times X\),
    \[
  \begin{aligned}
 &|u|^p\to Q^2, \, \text{in }L^1(\R^N)\cap L^\infty(\R^N),\\
  &|u|^{p-2}u\to Q \, \text{in }L^2(\R^N)\cap L^\infty(\R^N),\\
 &I_\alpha*|u|^p\to I_\alpha*Q^2\, \text{in }L^\infty(\R^N).
    \end{aligned}
 \]
 The last assertion follows from
   \begin{equation}\label{eq:quadratic-convolution}
    \|I_\alpha*g\|_{L^\infty(\R^N)} \le C\left( \|g\|_{L^1(\R^N)} + \|g\|_{L^\infty(\R^N)} \right).
  \end{equation}

We also claim that
\begin{equation}\label{eq:quadratic-multiplier}
  \left\| (I_\alpha*|u|^p)|u|^{p-2} - I_\alpha*Q^2 \right\|_{L^\infty(\R^N)} \to0.
  \end{equation}
 Indeed, given \(\varepsilon>0\), choose \(R>0\) such that
  \(\sup_{|x|>R}(I_\alpha*Q^2)(x)<\varepsilon\). The functions \(|u|^{p-2}\) are uniformly bounded for
  \((p,u)\) close to \((2,Q)\). Since \(\min_{\overline{B_R}}Q>0\) and \(u\to Q\) uniformly,
\[
  \left\||u|^{p-2}-1\right\|_{L^\infty(B_R)}\to0.
  \]
  Consequently,
    \[
  \begin{aligned}
 &\left\|(I_\alpha*|u|^p)|u|^{p-2} -I_\alpha*Q^2\right\|_{L^\infty(\R^N)}\\
  &\quad\le C\left\|I_\alpha*|u|^p-I_\alpha*Q^2\right\|_{L^\infty(\R^N)}
  +\|I_\alpha*Q^2\|_{L^\infty(\R^N)} \left\||u|^{p-2}-1\right\|_{L^\infty(B_R)}\\
&\qquad+C\sup_{|x|>R}(I_\alpha*Q^2)(x) \le o(1)+C\varepsilon.
  \end{aligned}
    \]
  Letting \((p,u)\to(2,Q)\) and then \(\varepsilon\downarrow0\) proves~\eqref{eq:quadratic-multiplier}.

Finally,
   \[
\begin{aligned}
   D_uF(p,u)h ={}& h - (-\Delta+1)^{-1} \Big[ (p-1)(I_\alpha*|u|^p)|u|^{p-2}h + p|u|^{p-2}u
     \bigl( I_\alpha*(|u|^{p-2}uh) \bigr) \Big].
  \end{aligned}
  \]
   We claim that
\begin{equation}\label{eq:quadratic-derivative}
   \|D_uF(p,u)-D_uF(2,Q)\|_{X\to X}\to0, \quad \text{as} \quad (p,u)\to(2,Q),
\end{equation}
 where
  \[
  D_uF(2,Q)h = h-(-\Delta+1)^{-1} \left[(I_\alpha*Q^2)h+2Q(I_\alpha*(Qh))\right].
  \]
 By~\eqref{eq:quadratic-multiplier} and \(p\to2\),
 \[
 \left\| (p-1)(I_\alpha*|u|^p)|u|^{p-2} - I_\alpha*Q^2 \right\|_{L^\infty(\R^N)} \to0.
\]
  Hence
 \begin{equation}
\label{eq:conv-2-infty}
  \left\| \left[ (p-1)(I_\alpha*|u|^p)|u|^{p-2} - I_\alpha*Q^2 \right]h \right\|_{L^2(\R^N)\cap L^\infty(\R^N)} =
  o(1)\|h\|_X.
   \end{equation}
Moreover, by a directly decomposition,
    \begin{equation}
  \label{eq:identify-decomposition}
    \begin{aligned}
  &p|u|^{p-2}u \bigl(I_\alpha*(|u|^{p-2}uh)\bigr) - 2Q\bigl(I_\alpha*(Qh)\bigr)\\
 &\quad= (p-2)|u|^{p-2}u \bigl(I_\alpha*(|u|^{p-2}uh)\bigr) + 2\bigl(|u|^{p-2}u-Q\bigr) \bigl(I_\alpha*(|u|^{p-2}uh)\bigr)\\
  &\qquad+ 2Q\,I_\alpha* \left( \bigl(|u|^{p-2}u-Q\bigr)h \right).
 \end{aligned}
    \end{equation}
   Since \(|u|^{p-2}u\to Q\) in \(L^2(\R^N)\cap L^\infty(\R^N)\), the family
   \(|u|^{p-2}u\) is uniformly bounded in \(L^2\cap L^\infty\).
   Thus~\eqref{eq:quadratic-convolution} gives
   \[
   \begin{aligned}
 \left\| I_\alpha*(|u|^{p-2}uh) \right\|_{L^\infty(\R^N)} &\le C\left( \||u|^{p-2}uh\|_{L^1(\R^N)} + \||u|^{p-2}uh\|_{L^\infty(\R^N)}
  \right) \le C\|h\|_X,
  \end{aligned}
 \]
   uniformly for \((p,u)\) close to \((2,Q)\), while
\[
  \begin{aligned}
  \left\| I_\alpha* \left[ \bigl(|u|^{p-2}u-Q\bigr)h \right] \right\|_{L^\infty(\R^N)} &\le C \left(
   \||u|^{p-2}u-Q\|_{L^2} + \||u|^{p-2}u-Q\|_{L^\infty(\R^N)} \right) \|h\|_X\\
 &= o(1)\|h\|_X.
\end{aligned}
\]
  It follows that
\begin{equation}
  \label{eq:conv-2-infty-2}
 \begin{aligned}
 & \left\| p|u|^{p-2}u \bigl(I_\alpha*(|u|^{p-2}uh)\bigr) - 2Q\bigl(I_\alpha*(Qh)\bigr)
   \right\|_{L^2\cap L^\infty} = o(1)\|h\|_X.
   \end{aligned}
\end{equation}
Thus by \eqref{eq:conv-2-infty}, \eqref{eq:identify-decomposition},
\eqref{eq:conv-2-infty-2}, and
  \[\|(-\Delta+1)^{-1}f\|_X \le \|f\|_{L^2(\R^N)}+\|f\|_{L^\infty(\R^N)},
  \]
   we obtain \eqref{eq:quadratic-derivative}.

 Write
\[
  D_uF(2,Q)=I-K, ~\text{where}~
  Kh=(-\Delta+1)^{-1}\left[(I_\alpha*Q^2)h+2Q(I_\alpha*(Qh))\right].
  \]
   To prove compactness of \(K\) on \(X\), let \((h_j)\) be bounded in
    \(X\). After passing to a subsequence, \(h_j\rightharpoonup h\) in \(H^1\) and
   \(h_j\to h\) in \(L^2_{\mathrm{loc}}\), with \(h\in X\). Since \(I_\alpha*Q^2\)
   is bounded and tends to zero at infinity, Lemma~\ref{lem:decay} and local
    compactness give
  \[
(I_\alpha*Q^2)(h_j-h) +2Q\bigl(I_\alpha*(Q(h_j-h))\bigr)\to0 \quad\text{in }L^2(\R^N).
    \]
  This sequence is bounded in \(L^\infty\) by \eqref{eq:quadratic-convolution}. Interpolation therefore gives
  convergence in \(L^r\) for every finite \(r>2\). Taking \(r>\max\{2,N/2\}\), the resolvent estimates and
  \(W^{2,r}(\R^N)\hookrightarrow L^\infty(\R^N)\) yield \(Kh_j\to Kh\) in \(X\). Moreover,
  \(D_uF(2,Q)h=0\) implies \(h\in H^2_{\mathrm{rad}}(\R^N)\) and \(L_Qh=0\), so
  Theorem~\ref{thm:index-one-nondegeneracy} gives \(h=0\). The Fredholm alternative now yields
  \begin{equation}\label{eq:quadratic-invertibility}
    D_uF(2,Q):X\to X \quad\text{is an isomorphism}.
   \end{equation}

 By~\eqref{eq:quadratic-derivative} and~\eqref{eq:quadratic-invertibility}, after reducing \(\delta>0\),
  the map
\[
 u\longmapsto u-[D_uF(2,Q)]^{-1}F(p,u)
  \]
    has Lipschitz constant at most \(1/2\) on the closed ball \(\{u\in
    X:\|u-Q\|_X\le\delta\}\), uniformly for \(p\ge2\) close to \(2\). Since
    \(F(p,Q)\to0\), choose \(0<\varepsilon<\delta\) such that
 \[
  \left\|[D_uF(2,Q)]^{-1}F(p,Q)\right\|_X \le\frac{\delta}{2}, \qquad 2\le p<2+\varepsilon.
\]
  The map therefore sends this ball into itself. The contraction mapping theorem and
  continuity of \(F\) give a unique continuous branch in the ball,
   \begin{equation}\label{eq:quadratic-branch}
  [2,2+\varepsilon)\ni p\longmapsto u_p\in X, \qquad F(p,u_p)=0,\qquad u_2=Q.
  \end{equation}
  In particular, \(u_p\to Q\) in \(X\) as \(p\downarrow2\).

   We next prove that this branch remains positive. Since \((I_\alpha*Q^2)(x)\to0\) as
     \(|x|\to\infty\), \eqref{eq:quadratic-multiplier} and \eqref{eq:quadratic-branch} allow us to
     choose \(R>0\) and decrease \(\varepsilon\) so that
   \[
 (I_\alpha*|u_p|^p)(x)|u_p(x)|^{p-2}\le\frac12 \qquad
\bigl(|x|\ge R,\ 2\le p<2+\varepsilon\bigr).
    \]
  Since \(\min_{\overline{B_R}}Q>0\) and \(u_p\to Q\) uniformly,
  we may further decrease \(\varepsilon\) so that \(u_p>0\) on
  \(B_R\). Testing \(-\Delta u_p+u_p=(I_\alpha*|u_p|^p)|u_p|^{p-2}u_p\)
  with \(-u_p^-\) gives
   \[
\begin{aligned}
  \|\nabla u_p^-\|_{L^2(\R^N)}^2+\|u_p^-\|_{L^2(\R^N)}^2 &= \int_{\R^N\setminus B_R}
  (I_\alpha*|u_p|^p)|u_p|^{p-2}(u_p^-)^2\,dx\le\frac12\|u_p^-\|_{L^2(\R^N)}^2.
  \end{aligned}
\]
  Hence \(u_p^-=0\). Since \(u_p\not\equiv0\), the strong maximum principle gives \(u_p>0\) in \(\R^N\).
    Thus the branch consists of centered positive radial solutions.

 Finally, suppose that \(Q\) and \(\widetilde Q\) are two centered positive radial solutions at
   \(p=2\). Applying \eqref{eq:quadratic-branch} separately at \(Q\) and at \(\widetilde Q\),
   we obtain \(\varepsilon_Q>0\), \(\varepsilon_{\widetilde Q}>0\), and positive centered radial solutions
   \(u_p\), \( \widetilde u_p\) for \(2\le p<2+\varepsilon_Q\), \(2\le p<2+\varepsilon_{\widetilde Q}\),
   respectively, such that
 \[
  u_2=Q, \qquad \widetilde u_2=\widetilde Q, \qquad u_p\to Q, \qquad \widetilde u_p\to\widetilde Q
  \qquad\text{in }X \quad\text{as }p\downarrow2.
 \]

  Set \(\varepsilon=\min\{\varepsilon_Q,\varepsilon_{\widetilde Q},\delta\}\). Then, for every
    \(2<p<2+\varepsilon\), both \(u_p\) and \(\widetilde u_p\) are centered positive radial solutions
    of~\eqref{eq:equation} with the same parameters \((\alpha,p)\). By
    Lemmas~\ref{lem:endpoint} and \ref{lem:proper-continuation}, \(u_p=\widetilde u_p\) for
    \(2<p<2+\varepsilon\). Choosing \(p_j\downarrow2\), we obtain
  \[
 Q = \lim_{j\to\infty}u_{p_j} = \lim_{j\to\infty}\widetilde u_{p_j} = \widetilde Q.
    \]
 Thus the centered positive radial solution at \(p=2\) is unique.
\end{proof}
  
\subsection{Proof of Theorem~\ref{thm:unique}}
 \begin{proof}[Proof of Theorem~\ref{thm:unique}]
 Existence follows from~\cite[Theorem~1]{MV13}.

 Suppose first that \(p>2\). By Proposition~\ref{lem:endpoint} and the amplitude
  normalization in~\eqref{eq:unnormalized}, equation~\eqref{eq:equation} has
  exactly one centered positive solution when \(\alpha\) is sufficiently close to \(N\).
  Lemma~\ref{lem:proper-continuation} shows that the number of centered positive
 solutions is constant throughout
   \[
   \max\{0,(N-2)p-N\}<\alpha<N.
    \]
   It is therefore equal to one throughout this interval. For \(p=2\), uniqueness
   of the centered positive solution follows from
   Lemma~\ref{lem:quadratic-uniqueness}.

    Finally, Theorem~\ref{thm:radial-symmetry} shows that every positive
    \(H^1(\R^N)\) solution becomes radial and strictly decreasing about
    the origin after a translation. The uniqueness of the centered positive
    solution therefore gives uniqueness up to translations.
   \end{proof}

\noindent\textbf{Conflict of interest:} Authors state no conflict of interest.

\noindent\textbf{Data Availability Statement:} Data sharing is not applicable to this article as no
datasets were generated or analysed during the current study.

\noindent{\bf Acknowledgement.}
G-D.~Li was supported by NSFC (No.12561019). 
 The authors acknowledge the use of AI tools. All mathematical arguments
and proofs in the final manuscript were checked and written by the authors.

\providecommand{\href}[2]{#2}
\providecommand{\arxiv}[1]{\href{http://arxiv.org/abs/#1}{arXiv:#1}}
\providecommand{\url}[1]{\texttt{#1}}
\providecommand{\urlprefix}{DOI }


\end{document}